\documentclass[12pt,leqno]{article}
\usepackage{amsfonts}
\usepackage{amsmath, amsthm, amsfonts, amssymb, color}
\usepackage{mathrsfs}
\usepackage{color}
\usepackage[notcite,notref]{showkeys}

\newtheorem{thm}{Theorem}[section]

\theoremstyle{definition}

\newcommand{\scr}[1]{\mathscr #1}
\definecolor{wco}{rgb}{0.5,0.2,0.3}

\numberwithin{equation}{section} \theoremstyle{remark}

\newcommand{\ua}{\uparrow}

\title{{\bf Donsker-Varadhan Large Deviations for Path-Distribution Dependent SPDEs }\footnote{Supported in
 part by  NNSFC (11771326, 11831014).} }
\author{
{\bf Panpan Ren$^{b,c)}$,    Feng-Yu Wang$^{a,b)}$  }\\
\footnotesize{$^{a)}$ Center for Applied Mathematics, Tianjin University, Tianjin 300072, China}\\
 \footnotesize{$^{b)}$ Department of Mathematics,
Swansea University, Bay Campus, SA1 8EN, United Kingdom}\\
\footnotesize{$^{c)}$ Mathematical Institute,Woodstock Road, OX2 6GG, University of Oxford}\\
\footnotesize{ 673788@swansea.ac.uk, Panpan.ren@maths.ox.ac.uk; wangfy@tju.edu.cn, F.-Y.Wang@swansea.ac.uk}}
\begin{document}
\allowdisplaybreaks
\def\R{\mathbb R}  \def\ff{\frac} \def\ss{\sqrt} \def\B{\mathbf
B}
\def\N{\mathbb N} \def\kk{\kappa} \def\m{{\bf m}}
\def\ee{\varepsilon}\def\ddd{D^*}
\def\dd{\delta} \def\DD{\Delta} \def\vv{\varepsilon} \def\rr{\rho}
\def\<{\langle} \def\>{\rangle}
  \def\nn{\nabla} \def\pp{\partial} \def\E{\mathbb E}
\def\d{\text{\rm{d}}} \def\bb{\beta} \def\aa{\alpha} \def\D{\scr D}
  \def\si{\sigma} \def\ess{\text{\rm{ess}}}\def\s{{\bf s}}
\def\beg{\begin} \def\beq{\begin{equation}}  \def\F{\scr F}
\def\Ric{\mathcal Ric} \def\Hess{\text{\rm{Hess}}}
\def\e{\text{\rm{e}}} \def\ua{\underline a} \def\OO{\Omega}  \def\oo{\omega}
 \def\tt{\tilde}\def\[{\lfloor} \def\]{\rfloor}
\def\cut{\text{\rm{cut}}} \def\P{\mathbb P} \def\ifn{I_n(f^{\bigotimes n})}
\def\C{\scr C}      \def\aaa{\mathbf{r}}     \def\r{r}
\def\gap{\text{\rm{gap}}} \def\prr{\pi_{{\bf m},\varrho}}  \def\r{\mathbf r}
\def\Z{\mathbb Z} \def\vrr{\varrho} \def\lll{\lambda}
\def\L{\scr L}\def\Tt{\tt} \def\TT{\tt}\def\II{\mathbb I}
\def\i{{\rm in}}\def\Sect{{\rm Sect}}  \def\H{\mathbb H}
\def\M{\mathbb M}\def\Q{\mathbb Q} \def\texto{\text{o}} \def\LL{\Lambda}
\def\Rank{{\rm Rank}} \def\B{\scr B} \def\i{{\rm i}} \def\HR{\hat{\R}^d}
\def\to{\rightarrow}\def\l{\ell}\def\iint{\int}\def\gg{\gamma}
\def\EE{\scr E} \def\W{\mathbb W}
\def\A{\scr A} \def\Lip{{\rm Lip}}\def\S{\mathbb S}
\def\BB{\scr B}\def\Ent{{\rm Ent}} \def\i{{\rm i}}\def\itparallel{{\it\parallel}}
\def\g{{\mathbf g}}\def\Sect{{\mathcal Sec}}\def\T{\mathcal T}\def\BB{{\bf B}}
\def\f{\mathbf f} \def\g{\mathbf g}\def\BL{{\bf L}}  \def\BG{{\mathbb G}}
\def\Bd{{D^E}} \def\BdP{D^E_\phi} \def\Bdd{{\bf \dd}} \def\Bs{{\bf s}} \def\GA{\scr A}
\def\Bg{{\bf g}}  \def\Bdd{{\bf d}} \def\supp{{\rm supp}}\def\div{{\rm div}}
\def\ddiv{{\rm div}}\def\osc{{\bf osc}}\def\1{{\bf 1}}\def\BD{\mathbb D}\def\GG{\Gamma}
\maketitle

\begin{abstract}  As an important tool characterizing the long time behavior of Markov processes, the Donsker-Varadhan LDP (large deviation principle)   does not directly  apply to distribution dependent SDEs/SPDEs since the solutions are non-Markovian. We establish this type LDP for
 several different models of  distribution dependent SDEs/SPDEs which may also with memories, by comparing the original equations with the corresponding distribution independent ones.   As preparations,  the existence, uniqueness and exponential convergence are also  investigated for path-distribution dependent SPDEs which should be interesting by themselves.
 \end{abstract} \noindent
 \vskip 2cm

 \section{Introduction}

The LDP (large deviation principle)  is a fundamental tool characterizing the asymptotic behaviour of probability measures  $\{\mu_\vv\}_{\vv>0}$ on a topological space $E$, see  \cite{DZ} and references within.
 Recall that  $\mu_\vv$ for small $\vv>0$ is said to
satisfy the LDP  with speed $\lambda(\vv)\to +\infty$ (as $\vv\to 0$) and rate
function $I: E\to [0,+\infty]$, if $I$ has compact level sets (i.e. $\{I\le r\}$ is compact  for $r\in {\mathbb R}^+$),
and for any Borel subset $A$ of $E$,
$$
 -\inf_{ A^o} I\le \liminf_{\vv\to 0}
   \frac 1{\lambda(\vv)} \log \mu_{\vv}(A)\le \limsup_{\vv\to 0}
   \frac 1{\lambda(\vv)} \log \mu_{\vv}(A)\le
  -\inf_{\bar A} I,
$$
where $A^o$ and $\bar A$ stand for    the
interior and the closure of $A$ in $E$ respectively.  The following two different type LDPs have been studied in the literature. 
  \paragraph{The Freidlin-Wentzell type small noise LDP \cite{FW}:} $\mu_\vv$ stands for the distribution of the solution to a dynamic system  perturbed by a noise with small intensity $\vv>0$, i.e. SDE (stochastic differential equation) with small noise. In this case, $E$ is the path space for the solutions of the SDE.
This type LDP describes, as $\vv\to 0$,  the convergence of stochastic systems to the corresponding deterministic system.
\paragraph{The Donsker-Varadhan type long time LDP \cite{DV}:}  $\mu_\vv$ stands for the distribution of $L_{\vv^{-1}},$  where
$$L_{t}:=\ff 1 t \int_0^{t} \dd_{X(s)}\d s,\ \ t>0$$  is the   empirical measure for a stochastic process $\{X(t)\}_{t\ge 0}.$ This type LDP describes the   behaviour of  $L_t$ as $t\to\infty$. In this case, $E$ is the set of all probability measures on the state space of the process,
on which  both the weak topology (induced by bounded continuous functions) and the $\tau$-topology (induced by bounded measurable functions) are considered in the literature.

\

In this paper, we study the Donsker-Varadhan LDP for path-distribution dependent SDEs (stochastic differential equations) on a separable Hilbert space $\H$. Inspired by  Kac's programme for Vlasov systems in  kinetic theory \cite{[18]}, McKean \cite{[32]} introduced distribution dependent SDEs. According to Sznitman \cite{[36]}, under  the global Lipschtiz condition,  these type SDEs can be derived as the  limit of mean-field  particle systems when  the number of particles tends to infinity.
  Therefore, distribution dependent SDEs are also called Mckean-Vlasov SDEs and  mean-field SDEs.

 In applications,  the distribution of a stochastic process can be regarded as a macro property, while the path of the  process up to a time $t$ stands for the history of the  system before this time. Since the evolution of  a stochastic system may depend on both the macro environment and the history, it is reasonable to investigate  path-distribution dependent SDEs. Moreover, because in many cases the configuration space for particle systems is infinite-dimensional, we consider path-distribution dependent  SDEs on Hilbert spaces, and in this case the SDEs are   called SPDEs (stochastic partial differential equations).

In recent years, distribution dependent SDEs   have been intensively investigated. Among many other papers in this field,  \cite{RST} established the Freidlin-Wentzell LDP for distribution dependent SDEs.
  However, up to our best knowledge, there is no any result on the Donsker-Varadhan LDP for this type SDEs. Since   the solution  is  non-Markovian, existing results on the Donsker-Varadhan LDP   derived  for Markov processes do  not apply. Indeed, the definition of the rate function (the Donsker-Varadhan level 2 entropy function)  depends on the Markov property of the process, for which the law  of the process starting at an initial distribution $\nu$ is given by
$$P^\nu=\int_EP^x\nu(\d x),$$ where   $ P^x $ is the law of the process starting at $x$, see Subsection 3.2 for details.

\paragraph{Main idea of the study.} To establish the Donsker-Varadhan type LDP for a distribution dependent SDE/SPDE, 
  we choose a reference SDE/SPDE whose solution is Markovian  so that existing results on  the Donsker-Varadhan LDP apply. 
  By comparing the original equation with the reference one in the sense of LDP,  see Lemma \ref{L1.1} below, we establish the Donsker-Varadhan LDP  for the distribution dependent SDE/SPDE. To this end, we will assume that the original equation has a unique invariant probability measure $\bar\mu$, and take the reference equation to be the original one with $\bar\mu$ replacing the distribution variable.

\paragraph{The framework.}  For a measurable space $(E,\scr B)$, let $\scr P(E)$ denote the set of all probability measures on $E$. For an $E$-valued random variable $\xi$ on a probability space, let $\L_\xi\in\scr P(E)$ be the distribution of $\xi$.

For a separable Hilbert space  $\H$, let $\mathbb L(\H)$ be the class  of all bounded linear operators on $\H$, which is equipped with the operator norm $\|\cdot\|$. We will also use the Hilbert-Schmidt norm $\|\cdot\|_{HS}$. In general, for two separable Hilbert spaces $\H_1$ and $\H_2$,
$\mathbb L(\H_1;\H_2)$ stands for the space of all bounded linear operators from $\H_1$ to $\H_2.$ 

For a fixed constant $r_0\ge 0$, let
$\mathscr{C}=C([-r_0,0];\mathbb{H})$ be the space of all continuous maps from $[-r_0,0]$ to $\mathbb H$ equipped with the uniform
norm $$\|\xi\|_\infty:=\sup_{-r_0\leq\theta\leq0}|\xi(\theta)|,\ \ \xi\in\C.$$
Then $\C$ is a Polish space, which refers to the history of a stochastic differential system on $\H$ with memory length $r_0$. 
When $r_0=0$, the path space  $\C$ degenerates to $\H$.
For any map $h(\cdot)\in C([-r_0,\infty);\mathbb{H})$ and $t\ge 0$, the corresponding segment  $h_t\in \C$ is defined by
 $$h_t(r)=h(t+r),\ \ \ r\in[-r_0,0].$$
Let $W(t)$ be the cylindrical Brownian motion on $\tt\H$ under a complete
filtration  probability space $(\OO,\F,\{\F_t\}_{t\ge 0}, \P)$; that
is,
$$W(t)=\sum_{i=1}^\infty B_i(t)\tt e_i,\ \ t\ge 0$$ for an orthonormal basis $\{\tt e_i\}_{i\ge 1}$ on $\tt\H$ and a sequence of independent one-dimensional
 Brownian motions $\{B_i(t): t\ge 0\}_{i\ge 1} $
on $(\OO,\F,\{\F_t\}_{t\ge 0}, \P)$,
where $\F_0$ is rich enough such that for any $\pi\in\scr P(\C\times\C)$ there exists a  $\C\times\C$-valued random variable $\xi$ on $(\OO,\F_0,\P)$ such that $\L_{\xi}=\pi$.  Consequently, for any $p>0$ and two probability measures 
$$\mu_1,\mu_2\in  \scr P_p(\C):=\big\{\mu\in \scr P(\C): \|\mu\|_p:= \mu(\|\cdot\|_\infty^p)^{\ff 1 {p\lor 1}}<\infty\big\},$$ 
there exist two $\F_0$-measurable random variables $\xi_1,\xi_2$ on $\C$ such that 
$$\big(\E[\|\xi_1-\xi_2\|_\infty^p]\big)^{\ff 1 p} = \W_p(\mu_1,\mu_2):=\inf_{\pi\in\C(\mu_1,\mu_2)} \bigg(\int_{\C\times\C} \|\xi-\eta\|_\infty^p\pi(\d\xi,\d\eta)\bigg)^{\ff 1 p}.$$ Note that $\scr P_p(\C)$ is a Polish  space under the $L^p$-Wasserstein distance $\{\W_p\}^{p\land 1}.$ 

Now, consider the following path-distribution dependent  SPDE on $\H$:
\begin{equation}\label{*1}
\d X(t)=\{AX(t)+b(X_t,\L_{X_t})\}\d t+\si(\L_{X_t})\d W(t),~~t\ge 0,
\end{equation}
where $(A,\D(A))$ is a negative definite self-adjoint operator on $\H$, $$b:\C\times\scr P(\C)\to\H,\ \ \si:  \scr P(\C)\to\mathbb L(\tt\H;\H)$$ are
measurable.

In Section 3, a more general equation \eqref{E'}   will be solved for initial value $X_0\in L^p(\OO\to\C, \F_0,\P)$ for some $p>0$; i.e. $\L_{X_0}$ is in the space

However, to establish the  Donsker-Varadhan LDP using the comparing method   proposed in Theorem \ref{TM} below, we have to assume that the noise term only depends on the distribution $\L_{X_t}$ rather than the solution $X_t$.

 Let $X_t^\nu$ denote the mild segment solution with initial distribution  $\nu\in \scr P(\C)$, which is a continuous adapted process on $\C$, see Definition 3.1 below for details.  We study the long time LDP for the empirical measure
$$L_t^\nu:= \ff 1 t \int_0^t \dd_{X_s^\nu}\d s,\ \ t>0.$$

 \beg{defn} \label{D2.1} Let $\scr P(\C)$ be equipped with the weak topology, let $\scr A\subset \scr P(\C)$, and let $J: \scr P(\C)\to [0,\infty]$ have compact level sets, i.e. $\{J\le r\}$ is compact in $\scr P(\C)$ for any $r>0$.
 \beg{enumerate} \item[(1)] $\{L_t^\nu\}_{\nu\in \scr A}$ is said to satisfy the upper bound uniform LDP with rate function $J$, denoted by $\{L_t^\nu\}_{\nu\in \scr A}\in LDP_u(J),$ if for any closed   $A\subset \scr P(\C),$
 $$\limsup_{t\to\infty} \ff 1 t\sup_{\nu\in \scr A} \log \P(L_t^\nu\in A)\le - \inf_AJ.$$
 \item[(2)] $\{L_t^\nu\}_{\nu\in \scr A}$ is said to satisfy the lower  bound uniform LDP with rate function $J$, denoted by $\{L_t^\nu\}_{\nu\in \scr A}\in LDP_l(J),$ if for any open  $A\subset \scr P(\C),$
 $$\liminf_{t\to\infty} \ff 1 t\inf_{\nu\in \scr A} \log \P(L_t^\nu\in A)\ge - \inf_AJ.$$
 \item[(3)] $\{L_t^\nu\}_{\nu\in \scr A}$ is said to satisfy the uniform LDP with rate function $J$, denoted by
 $\{L_t^\nu\}_{\nu\in \scr A}\in LDP(J),$ if $\{L_t^\nu\}_{\nu\in \scr A}\in LDP_u(J)$ and $\{L_t^\nu\}_{\nu\in \scr A}\in    LDP_l(J).$
\end{enumerate} \end{defn}

The remainder of the paper is orgnized as follows.  In Section 2, we state the main results of the paper and illustrate them by specific  examples. To prove these results, in Section 3 we investigate the existence and uniqueness for
 path-distribution dependent  SDEs/SPDEs, and recall some results on
the Donsker-Varadhan LDP for Markov processes.  Finally,  the proofs of main results are addressed in Section 4.

\section{Main results and Examples}

We investigate the long time LDP for \eqref{*1} in the following   three situations respectively:
 \beg{enumerate} \item[1)] $r_0=0$ and $\H$ is finite-dimensional;
 \item[2)]   $r_0=0$ and $\H$ is infinite-dimensional;   \item[3)] $r_0>0$ and $\si$ is constant.\end{enumerate}
 When $r_0>0$ and $\si$ is non-constant, there is a difficulty to apply our comparison argument.   We   leave this for the future study.
 
 To state our main results, we recall the Feller property, the strong Feller property and the irreducibility for a (sub-) Markov operator $P$.
 Let $\B_b(\C)$ (resp. $C_b(\C)$) be the space of bounded measurable (resp. continuous) real functions on $\C$.
Let $P$ be  a sub-Markov operator on $\B_b(\C)$, i.e. it is a positivity-preserving linear operator with $P1\le 1$.     $P$ is called strong Feller if $P\B_b(\C)\subset C_b(\C)$, is called Feller if $PC_b(\C)\subset C_b(\C)$, and is called $\mu$-irreducible for some $\mu\in \scr P(\C)$ if
$\mu(1_AP1_B)>0$ holds for any $A,B\in \B(\C)$ with $\mu(A)\mu(B)>0.$  

  \subsection{ Distribution dependent SDE on $\R^d$ }
Let $r_0=0,$  $\H=\R^d$ and $\tt\H=\R^m$ for some $d,m\in \mathbb N$.
  In this case, we combine the linear term $Ax$ with the drift term $b(x,\mu)$,  so that \eqref{*1} reduces to
  \beq\label{E0} \d X(t)= b(X(t),\L_{X(t)})\d t + \si(\L_{X(t)}) \d W(t),\end{equation}
  where $b: \R^d\times \scr P_2(\R^d)\to \R^d,$   $\si:  \scr P_2(\R^d)\to \R^d\otimes\R^m $ and $W(t)$ is the $m$-dimensional Brownian motion. We assume
  \beg{enumerate}\item[{\bf $(H_1)$}] $b$ is continuous, $\si$ is  bounded and   continuous such that
  $$2\<b(x,\mu) -b(y,\nu),x-y\>+ \|\si(\mu)-\si(\nu)\|_{HS}^2\le -\kk_1|x-y|^2 +\kk_2\W_2(\mu,\nu)^2$$ holds for some constants $\kk_1>\kk_2\ge 0$ and all $x,y\in\R^d, \mu,\nu\in \scr P_2(\R^d)$.
\end{enumerate}
Under {$(H_1)$}, for any $X(0)\in L^2(\OO\to\R^d,\F_0,\P)$,  the equation \eqref{E0} has a unique solution, see \cite[Theorem 2.1]{W18} or Theorem \ref{EXU} in a more general framework. We write $P_t^*\mu=\L_{X(t)}$ if $\L_{X(0)}=\mu$.  By \cite[Theorem 3.1(2)]{W18}, $P_t^*$ has a unique invariant probability measure $\bar\mu\in \scr P_2(\R^d)$ such that
\beq\label{EXPO}\W_2(P_t^*\nu, \bar\mu)^2\le\e^{-(\kk_1-\kk_2)t} \W_2(\nu,\bar\mu)^2,\ \ t\ge 0, \nu\in\scr P_2(\R^d).\end{equation}
Consider the reference SDE
\beq\label{E0'} \d \bar X(t)= b(\bar X(t),\bar\mu)\d t + \si(\bar\mu) \d W(t).\end{equation}
It is standard that under $(H_1)$   the equation \eqref{E0'} has a unique solution $\bar X^x(t)$ for any starting point $x\in \R^d,$ and $\bar\mu$ is the unique invariant probability measure of the associated Markov semigroup
$$\bar P_t f(x):= \E[f(\bar X^x(t))],\ \ t\ge 0, x\in \R^d, f\in \B_b(\R^d).$$
Consequently, $\bar P_t$ uniquely extends to $L^\infty(\bar\mu)$. If $f\in L^\infty(\bar\mu)$ satisfies
$$\bar P_t f= f+\int_0^t \bar P_s g\d s,\ \ \bar\mu\text{-a.e.}$$ for some $g\in L^\infty(\bar\mu)$ and  all $t\ge 0$, we write
$f\in \D(\bar{\scr A})$ and denote $\bar{\scr A} f=g$. Obviously, we have $\D(\bar{\scr A})\supset C_c^\infty(\R^d):=\{f\in C_b^\infty(\R^d): \nn f\text{\ has\ compact\ support}\}$ and
$$\bar{\scr A} f(x)= \ff 1 2 \sum_{i,j=1}^d \{\si\si^*\}_{ij}(\bar\mu)\pp_{i}\pp_j f(x)+ \sum_{i=1}^\infty b_i(x,\bar\mu) \pp_if(x),\ \ f\in C_c^\infty(\R^d).$$

According to Section 3,   the Donsker-Varadhan level 2 entropy function $J$ for the diffusion process generated by $\bar{\scr A}$ has compact level sets in $\scr P(\R^d)$ under the $\tau$ and weak topologies, and by \eqref{JN} below we have
$$ J(\nu)= \beg{cases} \sup\big\{\int_{\R^d} \ff{-\bar{\scr A} f}{f} \d\nu:\   1\le f\in \D(\bar{\scr A})\big\}, &\text{if}\ \nu\ll\mu,\\
 \infty, &\text{otherwise}.\end{cases} $$

\beg{thm} \label{T01} Assume {\bf $(H_1)$}. For any $r,R>0$, let
$\scr B_{r,R}=\big\{\nu\in \scr P(\R^d): \nu(\e^{|\cdot|^r})\le R\big\}.$
\beg{enumerate}
\item[$(1)$] We have $\{L^\nu_t\}_{\nu\in \scr B_{r,R}}\in LDP_u(J)$  for all $r,R>0$.
If $\bar P_t$ is strong Feller and $\bar\mu$-irreducible for some $t>0$, then $\{L^\nu_t\}_{\nu\in \scr B_{r,R}}\in LDP(J)$  for all $r,R>0$.
\item[$(2)$] If there exist constants $\vv,c_1,c_2>0$ such that
 \beq\label{ABC}  \<x, b(x,\nu)\>  \le c_1-c_2|x|^{2+\vv},\ \ x\in \R^d,\nu\in \scr P_2(\R^d),\end{equation} then $\{L^\nu_t\}_{\nu\in \scr P_2(\R^d)} \in LDP_u(J)$. If moreover $\bar P_t$ is strong Feller and $\bar\mu$-irreducible  for some $t>0$, then $\{L^\nu_t\}_{\nu\in \scr P_2(\R^d)}\in LDP(J).$   \end{enumerate}
\end{thm}

To apply this result, we first recall some facts on the strong Feller property and the irreducibility of diffusion semigroups. 
\paragraph{Remark 2.1.} (1) Let $\bar P_t$ be the (sub-)Markov semigroup generated by the second order differential operator 
$$\bar{\scr A}:= \sum_{i=1}^m U_i^2 +U_0,$$ 
where $\{U_i\}_{i=1}^m$ are $C^1$-vector fields and $U_0$ is a continuous vector field. According to   \cite[Theorem 5.1]{LAP}, if       $\{U_i: 1\le i\le m\}$    together with their Lie brackets with $U_0$  span $\R^d$ at any point (i.e. the H\"ormander condition holds),  then  the Harnack inequality
$$ P_t f(x)\le \psi(t,s,x,y)   P_{t+s} f(y),\ \ t,s>0, x,y\in \R^d, f\in \B^+(\R^d)$$
for some map $\psi: (0,\infty)^2\times (\R^d)^2\to (0,\infty).$ Consequently, if moreover $\bar P_t$ has an invariant probability measure $\bar\mu$, then  $\bar P_t$ is $\bar\mu$-irreducible for any $t>0.$   Finally, if $\{U_i\}_{0\le i\le m}$ are smooth with bounded derivatives of all orders, then  
the above H\"ormander condition implies that $\bar P_t$ has smooth heat kernel  with respect to the Lebesgue measure, in particular  it is strong Feller for any $t>0.$ 

(2) Let $\bar P_t$ be the  Markov semigroup generated by 
$$\bar{\scr A}:= \sum_{i,j=1}^d \bar a_{ij} \pp_i\pp_j  +\sum_{i=1}^d \bar b_i \pp_j,$$ 
where $(\bar a_{ij}(x))$ is strictly positive definite for any $x$, $\bar a_{ij}\in H_{loc}^{p,1}(\d x)$ and $\bar b_i\in L_{loc}^p(\d x)$ for some $p>d$ and 
all $1\le i,j\le d.$ Moreover, let $\bar\mu$ be an invariant probability measure of $\bar P_t$. Then by \cite[Theorem 4.1]{BKR},
$\bar P_t$ is strong Feller for all $t>0$. Moreover, as indicated in (1) that   \cite[Theorem 5.1]{LAP} ensures the $\bar\mu$-irreducibility of $\bar P_t$ for $t>0$.

  \
 
We present below two examples to illustrate this result, where the first is a distribution dependent perturbation of the Ornstein-Ulenbeck process, and the second is the distribution dependent stochastic Hamiltonian system.

\paragraph{Example 2.1.} Let $\si(\nu)=I+ \vv\si_0(\nu)$ and $b(x,\nu)= -\ff 1 2 (\si\si^*)(\nu)x$, where $I$ is the identity matrix, $\vv>0$ and $\si_0$ is a bounded Lipschitz continuous  map from $\scr P_2(\R^d)$ to $\R^d\otimes\R^d$. When $\vv>0$ is small enough, assumption  {\bf $(H_1)$} holds and 
that $\bar P_t$ satisfies conditions in Remark 2.1(2).  So, Theorem \ref{T01}(1)   implies    $\{L^\nu_t\}_{\nu\in \scr B_{r,R}}\in LDP(J)$  for all $r,R>0$,
where it is easy to see that the unique invariant probability measure $\bar\mu$ is the standard Gaussian measure on $\R^d$. 

If we take $b(x,\nu)= - x- c |x|^\theta x$ for some constants
$c,\theta>0$, then  when $\vv>0$ is small enough   $(H_1)$ and \eqref{ABC} are satisfied, so that Theorem \ref{T01}(2) and  Remark 2.1(2) imply $\{L^\nu_t\}_{\nu\in \scr P_2(\R^d)}\in LDP(J).$

\paragraph{Example 2.2.} Let $d=2m$ and consider the following distribution dependent SDE for $X(t)=(X^{(1)}(t), X^{(2)}(t))$ on $\R^{m}\times \R^m:$
$$\beg{cases} \d X^{(1)}(t)= \{X^{(2)}(t)- \lll X^{(1)}(t)\}\d t\\
\d X^{(2)}(t)= \{Z(X(t),\L_{X(t)}) - \lll X^{(2)}(t)\}\d t+\si\d W(t),\end{cases},$$
were $\lll>0$ is a constant, $\si$ is an invertible $m\times m$-matrix, $W(t)$ is the $m$-dimensional Brownian motion, and $Z:\R^{2m}\times\scr P_2(\R^{2m})\to\R^m$ satisfies
$$|Z(x_1,\nu_1)-Z(x_2,\nu_2)|\le \aa_1|x_1^{(1)}-x_2^{(1)}|+\aa_2|x_1^{(2)}-x_2^{(2)}|+\aa_3\W_2(\nu_1,\nu_2)$$ for some constants $\aa_1,\aa_2,\aa_3\ge 0$ and all
  $ x_i=(x_i^{(1)}, x_i^{(2)})\in \R^{2m}, \nu_i\in \scr P_2(\R^{2m}), 1\le i,j\le 2.$
  If
\beq\label{*PW} 4\lll> \inf_{s>0}\big\{2\aa_3 s + \aa_3s^{-1} + 2\aa_2+\ss{4(1+\aa_1)^2 + (2\aa_2+\aa_3s^{-1})^2}\big\},\end{equation}
then $\{L_t^\nu\}_{\nu\in \scr B_{r,R}} \in LDP(J)$ for all $r,R>0$.

Indeed,   $b(x,\nu):= (x^{(2)}-\lll x^{(1)}, Z(x,\nu)-\lll x^{(2)})$ satisfies
\beg{align*} & 2 \<b(x_1,\nu_1)-b(x_2,\nu_2), x_1-x_2\> \\
&\le -2\lll |x_1^{(1)}-x_2^{(1)}|^2 -2(\lll-\aa_2) |x_1^{(2)}-x_2^{(2)}|^2 \\
&\qquad + 2 |x_1^{(2)}-x_2^{(2)}|\big\{(1+\aa_1)|x_1^{(1)}-x_2^{(1)}| +\aa_3\W_2(\nu_1,\nu_2)\big\}\\
&\le \aa_3 s \W_2(\nu_1,\nu_2)^2 -\{2\lll-\dd(1+\aa_1)\}|x_1^{(1)}-x_2^{(1)}|^2\\
&\qquad-\{2\lll-2\aa_2-\dd^{-1}(1+\aa_1)-\aa_3s^{-1}\}|x_1^{(2)}-x_2^{(2)}|^2,\ \ s,\dd>0\end{align*}
for all $x_1,x_2\in\R^{2m}$ and $\nu_1,\nu_2\in \scr P_2(\R^{2m}).$
Taking
$$\dd= \ff{2\aa_2+\aa_3s^{-1} +\ss{4(1+\aa_1)^2+ (2\aa_2+\aa_3r^{-1})^2}}{2(1+\aa_1)}$$
such that $\dd(1+\aa_1)= 2\aa_2+\dd^{-1}(1+\aa_1)+\aa_3 s^{-1},$
we see that {\bf $(H_1)$} holds for some $\kk_1>\kk_2$ provided $2\lll-\dd(1+\aa_1)>\aa_3 s$ for some $s>0$, i.e. \eqref{*PW} implies {\bf $(H_1)$}.
Moreover,  it is easy to see that conditions in Remark 2.1(1) hold, see also \cite{GW12,WZ13} for  Harnack inequalities and gradeint estimates on stochastic Hamiltonian systems which also imply the strong Feller and $\bar\mu$-irreducibility of $\bar P_t$.      Therefore, the claimed assertion follows from Theorem  \ref{T01}(1).

\subsection{Distribution dependent SPDE}

Consider the following distribution-dependent SPDE on a separable Hilbert space $\H$:
  \beq\label{EP} \d X(t)= \{AX(t)+b(X(t),\L_{X(t)})\}\d t +\si(\L_{X(t)})\d W(t),\end{equation}
  where $(A,\D(A))$ is a linear operator on $\H$, $b: \H\times\scr P_2(\H)\to\H$ and $\si:  \scr P_2(\H)\to\mathbb L(\tt\H;\H)$ are measurable, and 
  $W(t)$ is the cylindrical Brwonian motion on $\tt\H$.
   We make the following assumption.

  \beg{enumerate} \item[{\bf $(H_2)$}]
  ($-A, \mathscr{D}(A))$ is    self-adjoint  with
  discrete spectrum
 $0<\lll_1\le \lll_2\le \cdots $ counting  multiplicities such that
 $\sum_{i=1}^\infty  \lll_i^{\gg-1}<\infty$ holds for some constant $\gg\in (0,1)$.

Moreover,  $b$ is Lipschitz  continuous on $\H\times\scr P_2(\H)$, $\si$ is bounded and there exist constants $\aa_1,\aa_2\ge 0$ with $\lll_1>\aa_1+\aa_2$ such that
$$2\<x-y, b(x,\mu)-b(y,\nu)\>+\|\si(\mu)-\si(\nu)\|_{HS}^2 \le 2\aa_1 |x-y |^2  +2\aa_2 \W_2(\mu,\nu)^2$$ holds for all $x,y\in \H$ and $ \mu,\nu\in \scr P_2(\H).$
\end{enumerate}
According to Theorem \ref{EXU} below, assumption {\bf $(H_2)$} implies
  that for any $X(0)\in L^2(\OO\to\H,\F_0,\P)$, the equation \eqref{EP} has a unique mild solution $X(t)$. As before we denote by $X^\nu(t)$ the solution with initial distribution $\nu\in \scr P_2(\H)$, and write $P_t^*\nu=\L_{X^\nu(t)}$.
Moreover, by It\^o's formula  and  $\kk:=\lll_1-(\aa_1+\aa_2)>0$, it is easy to see that $P_t^*$ has a unique invariant probability measure $\bar\mu\in \scr P_2(\H)$ and
\beq\label{EXPP} \W_2(P_t^*\nu,\bar\mu)\le \e^{-\kk t} \W_2(\nu,\bar\mu),\ \ t\ge 0.\end{equation}
Consider the reference SPDE
$$\d \bar X(t)= \{A\bar X(t)+b(\bar X(t),\bar\mu)\}\d t +\si(\bar\mu)\d W(t),$$ which is again well-posed
for any initial value $\bar X(0)\in L^2(\OO\to\H,\F_0,\P)$. Let $J$ be the Donsker-Varadhan level 2 entropy function for the Markov process $\bar X(t)$, see Section 3. For any $r,R>0$ let
$$\scr B_{r,R}:=\big\{\nu\in \scr P(\H): \nu(\e^{|\cdot|^r})\le R\big\}.$$

\beg{thm}\label{T02} Assume {\bf $(H_2)$}. If there exist constants $\vv\in (0,1)$ and $c>0$ such that
\beq\label{CCO}  \<(-A)^{\gg-1}x, b(x,\mu)\> \le c+ \vv|(-A)^{\ff \gg 2}x|^2,\ \  x\in\D((-A)^{\ff \gg 2}),  \end{equation}
then  $\{L^\nu_t\}_{\nu\in \scr B_{r,R}}\in LDP_u(J)$   for all $r,R>0$. If moreover $\bar P_t$ is strong Feller and $\bar\mu$-irreducible for some $t>0$, then $\{L^\nu_t\}_{\nu\in \scr B_{r,R}}\in LDP(J)$   for all $r,R>0$.
\end{thm}

Assumption {\bf $(H_2)$} is standard to imply  the well-posedness of \eqref{EP} and the exponential convergence of $P_t^*$ in $\W_2$.
 Condition \eqref{CCO} is implied by
\beq\label{GG} |(-A)^{\ff \gg 2 -1}b(x,\mu)|\le \vv' |(-A)^{\ff \gg 2}x|+c',\ \ x\in\D((-A)^{\ff \gg 2})\end{equation}  for some constants $\vv'\in (0,1)$ and $c'>0$. In particular,   \eqref{CCO} holds if $|b(x,\mu)|\le c_1+c_2|x|$ for some constants $c_1>0$ and $c_2\in (0, \lll_1).$
 When $\si$ is invertible with bounded $\si^{-1}$ and $b(\cdot,\mu)$ is Lipschitz continuous, 
   the dimension-free Harnack inequality established in \cite[Theorem 3.4.1]{W13} implies the strong Feller property and $\bar\mu$-irreducibility of
    $\bar P_t$ for $t>0$, see \cite[Theorem 1.4.1]{W13} for more properties implied by this type Harnack inequality.  Therefore, by Theorem \ref{T02},
   in this case $(H_2)$ and \eqref{GG} imply  $\{L^\nu_t\}_{\nu\in \scr B_{r,R}}\in LDP(J)$   for all $r,R>0$.
 See   Example 2.4 below for the case where $\si$ is non-invertible and $b$ is possibly  path-dependent.  

  \subsection{Path-distribution dependent SPDE with additive noise}

Let $\tt\H=\H$ and $\si\in \mathbb L(\H).$     Then \eqref{*1}  becomes
\beq\label{**0} \d X(t)= \big\{AX(t)+ b(X_t,\L_{X_t})\big\}\d t + \si \d W(t).\end{equation}
Below we consider this equation with  $\si$ being invertible and  non-invertible respectively.

\subsubsection{Invertible $\si$}
Since $\si$ is constant, we are able to establish LDP for $b(\xi,\cdot)$ being Lipshcitz continuous in $\W_p$ for some $p\ge 1$ rather than just for $p=2$ as in the last two results.

\begin{enumerate}
\item[{\bf $(H_3)$}] $\si\in\mathbb L(\H)$ is constant and $(A,\D(A))$ satisfies the corresponding condition in {\bf $(H_2)$}. Moreover, there exist constants $p\ge 1$ and  $\aa_1,\aa_2\ge 0$ such that
$$|b(\xi,\mu)-b(\eta,\nu)|\le \aa_1\|\xi-\eta\|_\infty +\aa_2 \W_p(\mu,\nu),\ \ \xi,\eta\in\C, \mu,\nu\in \scr P_p(\C).$$
 \end{enumerate}

Obviously,  {\bf $(H_3)$} implies assumption {\bf (A)} in Theorem \ref{EXU} below,  so that for any $X_0^\nu\in L^p(\OO\to\C,\F_0,\P)$ with $\nu=\L_{X_0^\nu}$,
the equation \eqref{*1} has a unique mild segment solution $X_t^\nu$ with
$$\E\Big[\sup_{t\in [0,T]} \|X_t^\nu\|_\infty^p\Big]<\infty,\ \ T>0.$$
 Let $P^*_t\nu=\L_{X_t^\nu}$ for $t\ge 0$ and $\nu\in \scr P_p(\C)$.

When $P_t^*$ has a unique invariant probability measure $\bar\mu\in \scr P_p(\C)$,  we consider the reference  functional SPDE
\beq\label{**} \d\bar X(t)= \big\{A\bar X(t) + b(\bar X_t, \bar\mu)\big\}\d t + \si\d W(t).\end{equation}
By Theorem \ref{EXU} below, this reference equation is well-posed for any initial value in $L^p(\OO\to \C,\F_0,\P)$.
For any $\vv,R>0$, let
$$\scr I_{\vv,R}=\big\{\nu\in \scr P(\C):  \nu(\e^{\vv\|\cdot\|_\infty^2}) \le R\big\}.$$

 \begin{thm}\label{TL1}
 Assume  {\bf $(H_3)$}.    Let $\theta\in [0,\lll_1]$ such that $$\kk_p:=
\theta-(\aa_1+\aa_2)\e^{p\theta r_0}=  \sup_{r\in [0,\lll_1]} \big\{r- (\aa_1+\aa_2)\e^{prr_0}\big\}.$$
 \beg{enumerate} \item[$(1)$] For any $\nu_1,\nu_2\in \scr P_p(\C)$,
 \beq\label{ES1} \W_p(P_t^*\nu_1, P_t^*\nu_2)^p\le \e^{p\theta r_0-p\kk_p t} \W_p(\nu_1,\nu_2)^p,\ \ t\ge 0.\end{equation} In particular, if $\kk_p>0$, then   $P_t^*$ has a unique invariant probability measure $\bar\mu\in \scr P_p(\C)$   such that
\beq\label{ES2} \W_p(P_t^*\nu, \bar\mu)^p\le \e^{p\theta r_0-p\kk_p t} \W_p(\nu,\bar\mu)^p,\ \ t\ge 0, \nu\in \scr P_p(\C).\end{equation}
 \item[$(2)$] Let $\si$ be invertible. If $\kk_p>0$ and $\sup_{s\in (0,\lll_1]}(s-\aa_1\e^{sr_0})>0$, then  $\{L_t^\nu\}_{\nu\in \scr I_{\vv,R}}\in LDP(J)$ for any $\vv,R>0$, where $J$ is the Donsker-Varadhan level 2 entropy function for the Markov process $\bar X_t$ on $\C$.
     \end{enumerate} \end{thm}

\paragraph{Example 2.3.} For a bounded domain $D\subset \R^d$, let $\H=L^2(D;\d x)$ and $A=-(-\DD)^\aa$, where $\DD$ is the Dirichlet  Laplacian on $D$ and
$\aa>\ff d 2$ is a constant. Let $\si=I$ be the identity operator on $\H$, and
$$b(\xi,\mu)= b_0(\mu)+ \aa_1 \int_{-r_0}^0 \xi(r)\Theta(\d r),\ \ (\xi,\mu)\in \C\times\scr P_1(\C),$$
where $\aa_1\ge 0$ is a constant,   $\Theta$ is a signed measure on  $[-r_0,0]$  with total variation $1$ (i.e. $|\Theta|([-r_0,0])=1$), and
$b_0$ satisfies
$$|b_0(\mu)-b_0(\nu)|\le \aa_2 \W_1(\mu,\nu),\ \ \mu,\nu\in\scr P_1(\C)$$ for some constant $\aa_2\ge 0$. Then {\bf $(H_3)$} holds for $p=1$, and as shown in he proof of Example 1.1 in \cite{BWY15} that $$\lll_1\ge \lll:=\ff{(d\pi^2)^{\aa}}{R(D)^{2\aa}},$$
where $R(D)$ is the diameter of $D$.
Therefore, all assertions in Theorem \ref{TL1} hold provided
$$\sup_{r\in (0,\lll]} \{r-(\aa_1+\aa_2)\e^{rr_0}\}>0.$$ In particular, under this condition   $\{L_t^\nu\}_{\nu\in \scr I_{\vv,R}}\in LDP(J)$  for any $\vv,R>1.$

\subsubsection{Non-invertible $\si$}

Let $\H=\H_1\times\H_2$ for two separable Hilbert spaces $\H_1$ and
$\H_2$, and consider the following path-distribution dependent SPDE for $X(t)=(X^{(1)}(t),X^{(2)}(t))$ on $\H$:
\beq\label{E1} \beg{cases} \d X^{(1)}(t)= \{A_1X^{(1)}(t)+  BX^{(2)}(t)\}\d t,\\
 \d X^{(2)} (t)= \{A_2X^{(2)}(t)+ Z(X_t,\L_{X_t})\}\d t +\si\d W(t),
\end{cases}
\end{equation} where  $(A_i,\D(A_i))$ is a densely defined closed linear operator  on $\H_i$  generating a $C_0$-semigroup $\e^{t A_i}$  ($i=1,2$), $B\in \mathbb L(\H_2; \H_1)$,
$Z: \C\mapsto \H_2$ is measurable, $\si\in\mathbb L(\H_2)$, and  $W(t)$ is the cylindrical
Wiener process on $\H_2$.  Obviously,   \eqref{E1} can be reduced to  \eqref{**0} by taking $A={\rm diag}\{A_1,A_2\}$ and using ${\rm diag}\{0,\si\}$ replacing $\si$,  i.e.  \eqref{E1}  is a special case of  \eqref{**0} with non-invertible $\si$.

For any $\aa>0$ and  $p\ge 1$,   define $$ \W_{p,\aa}(\nu_1,\nu_2):=\inf_{\pi\in\C(\nu_1,\nu_2)}
\bigg(\int_{\C\times\C} \big(\aa\|\xi_1^{(1)}-\xi_2^{(1)}\|_\infty+ \|\xi_1^{(2)}-\xi_2^{(2)}\|_\infty\big)^p \pi(\d\xi_1,\d\xi_2)\bigg)^{\ff 1 p}.$$
 We assume

\paragraph{$(H_4)$} Let $p\ge 1$ and $\aa>0$.
\beg{enumerate}
\item[$(H_4^1)$] ($-A_2, \mathscr{D}(A_2))$  is   self-adjoint  with
  discrete spectrum
 $0<\lll_1\le \lll_2\le \cdots $ counting  multiplicities such that $\sum_{i=1}^\infty \lll_i^{\gg-1}<\infty$ for some $\gg\in (0,1)$. Moreover, $ A_1 \le \dd-\lll_1 $   for some constant  $ \dd\ge 0$; i.e., $\<A_1 x,x\>\le (\dd-\lll_1)|x|^2$ holds for all $x\in \D(A_1)$.
\item[$(H_4^2)$] There exist constants $K_1,K_2> 0$ such that
\begin{equation*}
|Z(\xi_1,\nu_1)-b(\xi_2,\nu_2)|\le
K_1\|\xi_1^{(1)}-\xi_2^{(1)}\|_\infty +K_2 \|\xi_1^{(2)}-\xi_2^{(2)}\|_\infty+K_3 \W_{p,\aa}(\nu_1,\nu_2),~~(\xi_i,\nu_i) \in\C\times\scr P_p(\C).
\end{equation*}
\item[$(H_4^3)$] $\si$ is invertible on $\H_2$, and there exists $A_0\in \mathbb L(\H_1; \H_1)$ such that for any $t>0$, $B\e^{tA_2}= \e^{tA_1}\e^{tA_0} B$ holds   and
$$
Q_t:=\int_0^t\e^{sA_0}BB^*\e^{sA_0^*}\d s
$$
is invertible on $ \H_1$.
\end{enumerate}

By Theorem \ref{EXU} for $\H_0=\H_2$ and ${\rm diag}\{0,\si\}$ replacing $\si$,    $(H_4)$ implies that for any $X_0\in L^p(\OO\to\C,\F_0,\P)$ this equation has a unique mild segment solution. Let $P_t^*\nu=\L_{X_t}$ for $\L_{X_0}=\nu\in \scr P_p(\C).$

\begin{thm}\label{T03} Assume {\bf $(H_4)$} for some constants $p\ge 1$ and $\aa>0$ satisfying
\beq\label{K1}\aa\le \aa':=\ff 1 {2\|B\|} \big\{\dd-K_2+\ss{(\dd-K_2)^2+4K_1\|B\|}\big\},\end{equation}
where $\|\cdot\|$ is the operator norm. If
\beq\label{K2}   \inf_{s\in (0,\lll_1]}  s\e^{-  sr_0}>   K_2+\aa'\|B\| +  K_3,\end{equation}
then $P_t^*$ has a unique invariant probability measure $\bar\mu$ such that
\beq\label{K**} \W_p(P_t^*\nu,\bar\mu)^2\le  c_1\e^{-c_2 t} \W_p(\nu,\bar\mu),\ \ \nu\in \scr P_p(\C), t\ge 0 \end{equation}  holds for some constants $c_1,c_2>0$, and
  $\{L_t^\nu\}_{\nu\in \scr I_{\vv,R}}\in LDP(J)$  for any $\vv,R>1$, where $J$ is the Donsker-Varadhan level 2 entropy function for the associated reference equation for $\bar X(t)$.
\end{thm}

\paragraph{Example 2.4.} Consider the following equation for $X(t)= (X^{(1)}(t), X^{(2)}(t))$ on $\H=\H_0\times\H_0$ for a separable Hilbert space $\H_0$:
$$\beg{cases} \d X^{(1)}(t)= \{\aa_1 X^{(2)}(t)- \lll_1 X^{(1)}(t)\}\d t\\
\d X^{(2)}(t)= \{Z(X(t),\L_{X(t)}) - A   X^{(2)}(t)\}\d t+  \d W(t),\end{cases}$$
where $\aa_1\in\R\setminus\{0\}$, $W(t)$ is the cylindrical Brownian motion on $\H_0$,
$A$ is a self-adjoint operator on $\H_0$ with discrete spectrum such that all eigenvalues $0<\lll_1\le\ll_2\le\cdots $    counting multiplicities satisfy
$$\sum_{i=1}^\infty \lll_i^{\gg-1}<\infty$$ for some $\gg\in (0,1),$ and    $Z$ satisfies
$$|Z(\xi_1,\nu_1)-Z(\xi_2,\nu_2)|\le \aa_2\|\xi_1-\xi_2\|_\infty+\aa_3 \W_2(\nu_1,\nu_2),\ \ (\xi_i,\nu_i)\in \C\times\scr P_2(\C), i=1,2.$$
Let $$\aa = \ff 1 {2\aa_1} \Big(\ss{\aa_2^2+4\aa_1\aa_2}-\aa_2\Big).$$ Then $P_t^*$ has a unique invariant probability measure $\bar\mu\in\scr P_2(\C)$, and
  $\{L_t^\nu\}_{\nu\in \scr I_{R,q}}\in LDP(J)$  for any $R,q>1$ if
\beq\label{ASS} \inf_{s\in [0,\lll_1]} s\e^{-s r_0}> \aa_2+\aa_1\aa+\ff{\aa_3}{1\land\aa}.\end{equation}
  Indeed, it is easy to see that assumption {\bf $(H_4)$} holds for $p=2$, $\dd=0, \|B\|=\aa_1, K_1=K_2=\aa_2$ and $K_3=\ff{\aa_3}{1\land\aa}$. So, we have $\aa=\aa'$ and \eqref{ASS} is equivalent to \eqref{K2}. Then the desired assertion follows from
Theorem \ref{T03}.

 \section{Preparations }
In this part, we investigate path-distribution dependent SPDEs and  recall some facts on Donsker-Varadhan LDP for Markov processes.

\subsection{Path-distribution dependent SPDEs}
Consider the following path-distribution dependent SPDE on $\H$:
\beq\label{E'} \d X(t)= \big\{A X(t)+b_t(X_t,\L_{X_t})\big\}\d t+\si_t(X_t,\L_{X_t})\d W(t),\end{equation}
where  $(A,\D(A))$ is a negative self-adjoint operator  on $\H$,
and $$b: [0,\infty)\times \C\times \scr P(\C)\to \H,\ \ \si: [0,\infty) \times \C\times \scr P(\C)\to \mathbb L(\H;\tt\H)$$ are measurable, and $W(t)$ is the cylindrical Brownian motion on $\tt\H$.

\beg{defn}\label{DF} An   adapted continuous process $(X_t)_{t\ge 0}$ on $\C$ is called a mild segment (or functional) solution of \eqref{E'}, if
$$\E\int_0^t \big\{|\e^{(t-s)A} b_s(X_s,\L_{X_s})|+ \|\e^{(t-s)A} \si_s(X_s,\L_{X_s})\|_{HS}^2\big\}\d s  <\infty,\ \ t\ge 0,$$
and  the process $X(t):=X_t(0)$ satisfies $\P$-a.s.
$$X(t)= \e^{At}X(0)+ \int_0^t \e^{(t-s)A} b_s(X_s,\L_{X_s})\d s+\int_0^t \e^{(t-s)A}\si_s(X_s,\L_{X_s})\d W(s),\ \ t\ge 0.$$
In this case, we call $(X(t))_{t\ge 0}$ a mild solution of $\eqref{E'}$ with initial value $X_0$.
\end{defn}

To ensure the existence and uniqueness of mild solutions with $X_0\in L^p(\OO\to\C, \F_0,\P)$ for some $p>0$, we make the following assumption.

\paragraph{(A)}  Let $p\in (0,\infty)$. There exists a subspace $\H_0$ of $\H$ such that $\si(\xi,\nu)\tt\H\subset \H_0$ for any $(\xi,\nu)\in \C\times\scr P(\C),$ and the orthogonal projection $\pi_0: \H\to\H_0$ satisfies $A\pi_0=\pi_0A$  on $\D(A)$.
 Moreover, there exist  $\gg\in (0,1)$ and    $1\le K\in L^1_{loc} ([0,\infty)\to[0,\infty))$  such that
\beg{enumerate}\item[$(A_1)$]   $\int_0^ts^{-\gg} \|\e^{sA}\pi_0\|_{HS}^{2}\d s<\infty,\ t\in (0,\infty).$
\item[$(A_2)$]  There exists $p_0>2$ such that for any $t\ge 0, \xi,\eta\in\C$ and $ \mu,\nu\in \scr P_p(\C),$
\beg{align*}& |b_t(\xi,\mu)-b_t(\eta,\nu)|\le K(t) \big(\|\xi-\eta\|_\infty+\W_p(\mu,\nu)\big),\\
& \|\si_t(\xi,\mu)-\si_t(\eta,\nu)\|^{p}
 \le K(t)^{1\land \ff{p}{p_0}} \big(\|\xi-\eta\|_\infty^{p}+\W_p(\mu,\nu)^{p}\big).\end{align*}
 \item[$(A_3)$]  $ |b_t(0,\dd_0)|  + \|\si_t(0,\dd_0)\|^{p\lor p_0}   \le  K(t),\ \ t\ge 0.$
\end{enumerate}
In many references $(A_1)$ is replaced by $\int_0^ts^{-\gg} \|\e^{As}\|_{HS}^2\d s<\infty$, see for instance \cite{DP}. The present weaker version allows us to cover more examples with degenerate noise.

\paragraph{Remark 3.1.} By {\bf (A)} we have   $\e^{A(t-s)}\si_s=\e^{A(t-s)}\pi_0\si_s$, so that using $\e^{A(t-s)}\pi_0$ to replace the semigroup $S(t-s)$ in the proof of \cite[Proposition  7.9]{DP}, if $\Phi(s)$ is an adapted   process on $\mathbb L(\H;\tt\H)$
such that $\E \int_0^t\|\Phi(s)\|^q\d s<\infty$ for some $q>2$, then
$$W_\Phi(t):=  \int_0^t \e^{A(t-s)}\pi_0 \Phi(s) \d W(s),\ \ t\ge 0$$ is an adapted continuous  process on $\H$ such that
$$\E\bigg[\bigg|\inf_{s\in [0,t]}\int_0^t \e^{A(s-r)}\pi_0 \Phi(r) \d W(r)\bigg|^q\bigg]\le c \E \int_0^t\|\Phi(s)\|^q\d s$$
holds for some constant $c>0$.

\beg{thm}\label{EXU} Assume {\bf (A)} and let  $X_0\in L^p(\OO\to\C,\F_0,\P)$.
Then  $\eqref{E'}$ has a unique mild segment solution $\{X_t\}_{t\ge 0}$ starting at $X_0$ with
 $$  \E \Big[\sup_{t\in [0,T]} \|X_t\|^p_\infty\Big]<\infty,\ \ T\in (0,\infty), $$  provided one of  following conditions holds:
 \beg{enumerate} \item[$(1)$] $p>2$.
 \item[$(2)$] $p\in (0,2]$ and $\si_s(\xi,\mu)$ does not depend on $\xi$.
 \item[$(3)$] $p=2$ and for any $s\ge 0, \xi,\eta\in \C$ and $\mu,\nu\in \scr P_2(\C)$, $$\|\si_s(\xi,\mu)-\si_s(\eta,\nu)\|_{HS}^2\le K(s) \{\|\xi-\eta\|_\infty^2 +\W_2(\mu,\nu)^2\big\}.$$
     \end{enumerate}
\end{thm}

\beg{proof}  We consider cases (1)-(3) respectively. 

\ \newline 
\emph{Proof for Case (1).} Let $p>2$.
\newline 
{\bf The existence.}   We adopt an iteration argument as in \cite{W18}.
It suffices to prove that for any fixed $T>0$, the SPDE has a unique mild segment solution
up to time $T$ satisfying
\beq\label{ME}  \E \Big[\sup_{t\in [0,T]} \|X_t\|^p_\infty\Big]<\infty.\end{equation}

(1a) We first consider the case that   $X_0$ is bounded. Let $X_t^0=X_0$  and $\mu_t^0=\L_{X_t^0}$ for $t\ge 0$.
By Remark 3.1,  
$$X^1(t):= \e^{At}X(0)+\int_0^t \e^{(t-s)A} b_s(X_s^0,\mu_s^0)\d s+\int_0^t \e^{(t-s)A}\si_s(X_s^0,\mu_s^0)\d W(s),\ \ t\ge 0,$$ is an  adapted continuous  process on $\H$
such that
\beq\label{OOP}  \E \Big[\sup_{t\in [0,T]} \|X_t^1\|^q_\infty\Big]<\infty,\ \  q>0,\end{equation} where
$X_t^1(r):= X^1(t+r)1_{\{t+r\ge 0\}} + X_0(t+r)1_{\{t+r<0\}}.$

Now, assume that for some $n\ge 1$ we have constructed a continuous  adapted process $\{X^n_t\}_{t\in [0,T]}$ on $\C$ with $X_0^n=X_0$ and
$$ \E \Big[\sup_{t\in [0,T]} \|X_t^n\|^q_\infty\Big]<\infty,\ \ q>0.$$ By Remark 3.1, 
\beq\label{*P} X^{n+1}(t):=  \e^{At}X(0)+\int_0^t \e^{(t-s)A} b_s(X_s^n,\mu_s^n)\d s+\int_0^t \e^{(t-s)A}\si_s(X_s^n,\mu_s^n)\d W(s),\ \ t\in [0,T]\end{equation}
for $\mu_s^n:=\L_{X_s^n}$ is an adapted  continuous   process on $\H$, and   the segment process $X_t^{n+1}$ given by
\beq\label{*P'} X_t^{n+1}(r):= X^{n+1}(t+r)1_{\{t+r\ge 0\}} + X_0(t+r)1_{\{t+r<0\}}\ \text{for}\ r\in [-r_0,0], \ \ t\ge 0\end{equation} satisfies
$$ \E \Big[\sup_{t\in [0,T]} \|X_t^{n+1}\|^q_\infty\Big]<\infty,\ \ q>0.$$
It suffices to find a constant $t_0>0$ independent of $X_0$ such that   $\{X_{[0,t_0]}^n\}_{n\ge 1}$ is a Cauchy sequence in $L^p(\OO\to C([0,t_0];\C),\P)$. This together with  assumption  {\bf (A)}
imply that the limit $X_{[0,t_0]}:=\lim_{n\to\infty} X_{[0,t_0]}^n$  gives rise to a mild segment  solution of \eqref{E'} up to time $t_0$.
By repeating the procedure with initial time $it_0$ and initial value $X_{it_0}$ for $i\ge 1$, in finite many steps we may construct a mild segment solution of \eqref{E'} up to time $T$, such that \eqref{ME} holds.

For any $n\ge 1, $ by \eqref{*P}, \eqref{*P'} and assumption {\bf (A)} we have
\beq\label{*RR} \beg{split} &\psi_n(t):= \sup_{s\in [0,t]} \|X_s^{n+1}-X_s^n\|_\infty =\sup_{s\in [0,t]} |X^{n+1}(s)-X^n(s)| \\
&\le \int_0^t K(s)\big\{\|X_s^n-X_s^{n-1}\|_\infty + \W_p(\mu_s^n,\mu_s^{n-1})\big\}\d s + \sup_{s\in [0,t]}
\bigg|\int_0^s \e^{A(s-r)}\Phi_n(r) \d W(r)\bigg|,\end{split}\end{equation}
where $\Phi_n(r):= \si_s(X_r^n,\mu_r^n)- \si_s(X_r^{n-1},\mu_r^{n-1})$ satisfies
$$\|\Phi_n(r)\|^{p}\le K(r)\big\{\|X_r^n-X_r^{n-1}\|_\infty^{p} + \W_p(\mu_r^n,\mu_r^{n-1})^{p}\big\}.$$
Combining this with
$ \W_p(\mu_r^n,\mu_r^{n-1})^p\le \E \|X_r^n-X_r^{n-1}\|_\infty^p,$  and noting that
Remark 3.1 implies
\beq\label{*RR2} \E\bigg[\sup_{s\in [0,t]}
\bigg|\int_0^s \e^{A(s-r)}\Phi_n(r) \d W(r)\bigg|^{p} \bigg]\le c \E \int_0^t \|\Phi_n(s)\|^{p}\d s\end{equation} for some constant $c>0$, we find   constants $C_1,C_2>0$ such that
\beg{align*} \E[\psi_n(t)^p] &\le C_1 \E\bigg(\int_0^t K(s)  \big\{\psi_{n-1}(s)
 +\W_p(\mu_s^n,\mu_s^{n-1})\d s\bigg)^p\\
 &\quad+  C_1 \E\int_0^t K(s)\big\{\|X_s^n-X_s^{n-1}\|_\infty^{p} + \W_p(\mu_s^n,\mu_s^{n-1})^{p}\big\}\d s  \\
& \le C_2\vv(t) \E[\psi_{n-1}^p(t)],\ \ \vv(t):= \bigg(\int_0^tK(s)\d s\bigg)^p + \int_0^tK(s)\d s.\end{align*}
Taking $t_0\in (0,T]$ such that $C_2\vv(t)\le \ff 1 2$, we obtain
$$\E [\psi_n^p(t_0)]\le   2^{-n} \E[\psi_0^p(t_0)]<\infty,\ \ n\ge 1.$$
Thus, $\{X_{[0,t_0]}^n\}_{n\ge 1}$ is a Cauchy sequence in $L^p(\OO\to C([0,t_0];\C),\P)$ as desired.

(1b) In general, for $X_0\in L^p(\OO\to\C,\F_0,\P)$ and $N\in\mathbb N$ let $X_0^{(N)}= X_01_{\{\|X_0\|_\infty\le N\}}$. By (1a), for any $N\ge 1$ we have constructed a mild segment solution $(X_t^{(N)} )_{t\in [0,T]}$ for \eqref{E'} satisfying \eqref{ME} with initial value $X_0^{(N)}$:
$$X^{N}(t):=  \e^{At}X^{(N)}(0)+\int_0^t \e^{(t-s)A} b_s(X_s^{(N)},\mu_s^{(N)})\d s+\int_0^t \e^{(t-s)A}\si_s(X_s^{(N)},\mu_s^{(N)})\d W(s),\ \ t\in [0,T],$$
where $\mu_s^{(N)}= \L_{X_s^{(N)}}.$ By the above argument
for $X^{(N)}(t)-X^{(M)}(t)$ instead of $X^{n+1}(t)-X^n(t)$,  we find a constant $C>0$ such that for any $N,M\ge 1$, the process
 $$\psi_{N,M} (t):=\sup_{s\in [0,t]} \|X_s^{(N)}-X_s^{(M)}\|_\infty^p,\ \ t\in [0,T]$$ satisfies
 \beq\label{*PP} \E[\psi_{N,M}(t)] \le C \E\big[\|X_0\|_\infty^p1_{\{\|X_0\|_\infty>N\land M\}}\big]+C\vv(t)  \E[\psi_{N,M}(t)],\ \ t\in [0,T].\end{equation}
 Taking $t_0\in (0,T]$ such that $C\vv(t_0)\le \ff 1 2$, we obtain
 $$\E[\psi_{N,M}(t_0)] \le 2C \E\big[\|X_0\|_\infty^p1_{\{\|X_0\|_\infty>N\land M\}}\big],\ \ N,M\ge 1,$$
 so that,  $\{X_{[0,t_0]}^{(N)}\}_{N\ge 1}$ is a Cauchy sequence in $L^p(\OO\to C([0,t_0];\C),\P)$, and it is easy to see that its limit as $N\to\infty$ is a solution of \eqref{E'} up to time $t_0$. As explained before that by repeating the procedure we construct a mild segment solution of \eqref{E'} up to time $T$ satisfying \eqref{ME}.
\ \newline 
{\bf The uniqueness.} Let $X(t)$ and $Y(t)$ be two mild segment solutions with initial value $X_0$ satisfying
$$\E\bigg[\sup_{t\in [0,T]} \big(\|X_t\|_\infty^p+ \|Y_t\|_\infty^p\big)\bigg]<\infty.$$
Similarly to \eqref{*PP} we have
$$   \E\Big[\sup_{s\in [0,t]}\|X_s-Y_s\|_\infty^p \Big] \le   C\vv(t)  \E\Big[\sup_{s\in [0,t]}\|X_s-Y_s\|_\infty^p \Big],\ \ t\in [0,T].$$ This implies $X_t=Y_t$ up to time $t_0\in (0,T]$ such that $C\vv(t_0)<1.$   Since this $t_0$ does not depend on the initial value, repeating the same argument leads to $X_t=Y_t$ for all $t\in [0,T].$ 

\ \newline 
\emph{ Proof for  Case (2).}  Let $p\in (0,2]$. Again we first assume that $X_0$ is bounded and let $X^n, \mu^n,\psi_n$ be defined in step (1a). Since $\si_s(\xi,\mu)$ does not depend on $\xi$ and $K(s)\ge 1$, by $(A_2)$,
$\Phi_n(s)$ in \eqref{*RR} satisfies
$$\|\Phi_n(s)\|^{p_0}\le K(s)  \W_p(\mu_s^n,\mu_s^{n-1})^{p_0}.$$
Combining this with Remark 3.1 for $q=p_0>2$, and using $\W_p(\mu_s^n,\mu_s^{n-1})^p\le \E \|X_s^n-X_s^{n-1}\|_\infty^p$,   we find   a constant  $C_1>0$ such that
\beg{align*} \E[\psi_n^p(t)] &\le C_1 \bigg(\E\int_0^t K(s)  \psi_{n-1} (s) \d s\bigg)^p +  C_1 \bigg(\int_0^t K(s) \W_p(\mu_s^n,\mu_s^{n-1})^{p_0} \d s\bigg)^{\ff p{p_0}} \\
& \le C_1\dd(t) \E[\psi_{n-1}^p(t)],\ \ t\in [0,T], n\ge 1\end{align*} holds for
$\dd(t):= (\int_0^tK(s)\d s)^p+ (\int_0^tK(s)\d s)^{\ff p {p_0}}.$
Then the remainder of the proof, including the existence and uniqueness for bounded $X_0$, and the extension to general $X_0\in L^p(\OO\to\C,\F_0,\P)$,     is similar to that in Case (1).

\ \newline 
\emph{Proof for  Case (3).} Let $p=2$. As explained above we only consider bounded   $X_0$. In this case, let $X^n, \mu^n,\psi_n$ be defined in step (1a).
By {\bf (A)} and It\^o's formula to $|X^{n+1}(t)- X^n(t)|^2$,  we find a constant $c>0$ such that
\beq\label{LOOP} \beg{split}  &\d |X^{n+1}(t)-X^n(t)|^2\\
&\le c K(t) \Big\{  |X^{n+1}(t)-X^n(t)|^2 +   \|X^n_t- X^{n-1}_t\|^2_\infty+ \W_2(\mu_t^n,\mu_t^{n-1})^2\Big\}\d t +\d M^n(t),\end{split}\end{equation}
where   $$M^n(t):=2\int_0^t
   \big\<X^{n+1}(s)-X^n(s), \{\si_s(X^n_s,\mu_s^n)-\si_s(X_s^{n-1},\mu_s^{n-1})\}\d W(s)\big\>$$ satisfies
$$\d\<M^n(t)\>\le 4 K(t)  |X^{n+1}(t)-X^n(t)|^2\big\{\|X^n_t- X^{n-1}_t\|_\infty^2+ \W_2(\mu_t^n,\mu_t^{n-1})^2\big\}\d t.$$
Obviously, \eqref{LOOP} implies
\beg{align*} &|X^{n+1}(t)-X^n(t)|^2\e^{-c \int_0^tK(s) d s} \\
&\le    \int_0^tK(s) \e^{-c\int_0^sK(r)\d r} \big(\|X^n_s- X^{n-1}_s\|^2_\infty+ \W_2(\mu_s^n,\mu_s^{n-1})^2\big)\d s + \int_0^t \e^{-\int_0^s c K(r)\d r } \d M^n(s) \end{align*} for $t\in [0,T].$
Therefore, by the BDG inequality, there exist constants $C_1, C_2>0$ depending only on $ T$  such that
\beg{align*}& \E [\psi_n^2(t)]\le \e^{c\int_0^T K(s)\d s} \E\bigg[\sup_{s\in [0,t]}|X^{n+1}(s)-X^n(s)|^2\e^{-c \int_0^sK(r) d r}\bigg] \\
&\le   C_1  \int_0^t K(s) \big\{\E [\psi_{n-1}^2(s)]+   \W_2(\mu_s^n,\mu_s^{n-1})^2\big\}\d s \\
& \qquad +C_1 \E\bigg[\bigg(\int_0^t K(s) \psi_n^2(s)\big\{  \psi_{n-1}^2(s) +   \W_2(\mu_s^n,\mu_s^{n-1})^2\big\}\d s\bigg)^{\ff 1 2}\bigg]  \\
&\le \ff 1 2 \E [\psi_n^2(t)] + C_2 \int_0^t K(s) \big\{\E [\psi_{n-1}^2(s)]+  \W_2(\mu_s^n,\mu_s^{n-1})^2\big\}\d s,\ \ t\in [0,T].\end{align*}
Noting that $ \W_2(\mu_s^n,\mu_s^{n-1})^2\le \E[\psi_{n-1}^2(s)],$ this implies
$$\E [\psi_n^2(t)]\le 4C_2\int_0^tK(s)  \E[\psi_{n-1}^2(s)]\d s\le 4 C_2  \E[\psi_{n-1}^2(t)] \int_0^tK(s)\d s,\ \ t\in [0,T], n\ge 1.$$
Taking $t_0\in (0,T]$ such that $4C_2\int_0^{t_0}K(s)\d s\le \ff 1 2$, we obtain
$$\E [\psi_n^2(t_0)] \le 2^{-(n-1)} \E[\psi_0^2(t_0)]<\infty,\ \ n\ge 1.$$
Thus, $\{X_{[0,t_0]}^n\}_{n\ge 1}$ is a Cauchy sequence in $L^2(\OO\to C([0,t_0];\C),\P)$ as desired.
The remainder of the proof is similar to that in Case (1).
  \end{proof}

\subsection{Donsker-Varadhan LDP for Markov processes}

We first introduce the rate function, i.e. the Donsker-Varadhan level 2 entropy function for continuous  Markov processes on a Polish space $E$.

Consider the path space
 $${\bf C}_E:=C([0,\infty)\to E)=\{w: [0,\infty)\ni t\mapsto w(t)\in E \ \text{is\ continuous}\}.$$  Let $\scr P({\bf C}_E)$ be the set of all probability measures on ${\bf C}_E$,
and $\scr P^s({\bf C}_E)$  the set of all stationary (i.e. time-shift-invariant) elements in
$\scr P({\bf C}_E)$. For any $Q\in \scr P^s({\bf C}_E)$,  let $\bar Q$ be the unique stationary probability measure on $\bar {\bf C}_E:=C(\R\to E)$ such that
$$\bar Q\big(\{w\in \bar {\bf C}_E: w(t_i)\in A_i, 1\le i\le n\}\big)= Q\big(\{w\in {\bf C}_E: w(t_i+s)\in A_i,1\le i\le n\}\big)$$ holds for any $n\ge 1,  -\infty<t_1<t_2<\cdots<t_n<\infty, s\ge -t_1, $ and $\{A_i\}_{1\le i\le n}\subset   \B(E).$
We call $\bar Q$ the stationary extension of $Q$ to $\bar{\bf C}_E.$ For any $s\le t$, let $\F_t^s:=\si(\bar {\bf C}_E\ni w\mapsto w(u): s\le
u\le t)$. For a probability measure $\bar Q$ on $\bar {\bf C}_E$, let  $\bar Q_{w-}$ be the regular
conditional distribution of $\bar Q$ given $\F_0^{-\infty}$.  Moreover, let
$\Ent_{\F_1^0}$ be the   Kullback-Leibler divergence (i.e. relative entropy) on the
$\si$-field $\scr F_1^0$; that is, for any two probability measures $\mu_1,\mu_2$ on ${\bf C}_E$,
$$\Ent_{\scr F_1^0}(\mu_1|\mu_2):=\beg{cases}  \int_{ {\bf C}_E} \big(h \log h  \big)\d\mu_2,\ &\text{if}\ \d \mu_1|_{\F_1^0}= h \d \mu_2|_{\scr F_1^0},\\
\infty, &\text{otherwise}.\end{cases}$$

Now, for a standard Markov process on $E$ with   $\{P^x:x\in E\} \subset \scr P({\bf C}_E),$ where $P^x$ stands for the distribution of the   process starting at $x$,   the process level entropy function of
Donsker-Varadhan is given by
$$H (Q):=\beg{cases}  \int_{\bar {\bf C}_E}\Ent_{\F_1^0}(\bar Q_{w-}|
P^{w(0)}) \bar Q(\d w), &\text{if}\ Q\in \scr P^s({\bf C}_E),\\
\infty, &\text{otherwise.}\end{cases}$$    Then the Donsker-Varadhan level 2 entropy
function  is defined as
\beq\label{JJ} J(\nu):=\inf\big\{H(Q):\ Q\in \scr P^s({\bf C}_E), Q(w(0)\in\cdot)=\nu\big\},\ \ \
\nu\in\scr P (E).\end{equation} This function has compact level sets in $\scr P(E)$ under the $\tau$- (hence the weak) topology, see for instance
\cite{Wu00, Wu00b}.
For any $\nu\in \scr P(E)$, let   $(X_t^\nu)_{t\ge 0}$ be the Markov process with initial distribution $\nu$. Consider its empirical measure
$$L_t^\nu:=\ff 1 t \int_0^t \dd_{X_s^\nu}\d s,\ \ t>0.$$ When $\nu=\dd_x$, we denote $X_t^\nu=X_t^x$ and $L_t^\nu=L_t^x$.
  Let $\mu$ be an invariant probability measure of $P_t$, where $P_t$ is the Markov semigroup given by
$$P_t f(x)=\E[f(X_t^x)],\ \ x\in E, t\ge 0, f\in \B_b(E). $$
We write $f\in \D_\mu(\scr A)$ if $f\in L^\infty(\mu)$ and
there exists $g\in L^\infty(\mu)$ such that $P_tf-f=\int_0^t P_sg\d s$ holds $\mu$-a.e. for all $t\ge 0$.  In this case, we denote $\scr A f=g$.
We have the following formula for $J$.

\beg{thm}[\cite{Wu00b}, Proposition B.10 and Corollary B.11] \label{T0} Assume that $P_t$ has a unique invariant probability measure $\mu$. Then
\beq\label{JN} J(\nu)= \beg{cases} \sup\big\{\int_E \ff{-\scr A f}{f} \d\nu:\   1\le f\in \D_\mu(\scr A)\big\}, &\text{if}\ \nu\ll\mu,\\
 \infty, &\text{otherwise}.\end{cases} \end{equation}
In particular, if the Markov process is associated with a symmetric Dirichlet form $(\EE,\D(\EE))$ in $L^2(\mu)$, then
\beq\label{JN2}J(\nu)= \beg{cases} \EE(h^{\ff 1 2}, h^{\ff 1 2}), &\text{if}\ \nu=h\mu, h^{\ff 1 2}\in \D(\EE),\\
\infty, &\text{otherwise}.\end{cases} \end{equation} \end{thm}

We now recall another result due to \cite{Wu00b} on the LDP for uniformly integrable Markov semigroups, which will be used in the proof of
Theorem \ref{TL1}.   Let $p\ge 1$ and let $P$ be a bounded linear operator on $L^p(\mu)$. We call $P$ uniformly integrable  in $L^p(\mu)$ if
$$\lim_{R\to\infty} \sup_{\mu(|f|^p)\le 1} \mu(|P f|^p1_{\{|Pf|>R\}})=0.$$ This LDP is established under the $\tau$-topology induced by $f\in \B_b(E)$, and hence also holds under the weak topology. Let $\nu\in I_{q, L}:= \{\nu=h\mu: \|h\|_{L^q(\mu)}\le L\}$ for   $q,L\in (1,\infty).$

\beg{thm}[\cite {Wu00b}, Theorem 5.1] \label{T2} Assume that the Markov semigroup $P_t$ has a unique invariant probability measure $\mu$,  and there exists $T\in (1,\infty)$ and   $p\in (1,\infty)$ such that $P_T$ is $\mu$-irreducible and uniformly integrable in $L^p(\mu)$. Then $\{L_t^\nu\}_{\nu\in I_{q, L}}\in LDP(J)$   under the $\tau$-topology  for all $q,L\in (1,\infty).$  \end{thm}

The next result due to \cite{Wu00}  provides criteria on the LDP using the hitting time to compact sets, which  will be used in the proofs of Theorem \ref{T01} and Theorem \ref{T02}.  For any set $K\subset E$ and any $x\in E$, let
$$\tau_K^x:=\inf \{t\ge 0: X^x(t)\in K\},$$ where $X^x(t)$ is the Markov process starting at $x$. We will use the following conditions:
\beg{enumerate} \item[{\bf (D1)}] For any $\lll>0$ there exists a compact set $K\subset E$ such that
\beq\label{Wu1} \sup_{x\in E} \E [\e^{\lll \tau_K^x}]<\infty.\end{equation}
\item[{\bf (D2)}] For any $\lll>0$ there exist a constant  $ s>0$ and a compact set $K\subset E$ such that for any compact set $K'\subset E$,
\beq\label{Wu2} \sup_{x\in K} \E[\e^{\lll \tau_{K}^{X^x(s)}}]<\infty,\ \ \sup_{x\in K'} \E[\e^{\lll \tau_K^x}]<\infty.\end{equation}
\end{enumerate}

\beg{thm} [\cite{Wu00}, Theorems 1.1,1.2] \label{T3} Assume that $P_t$ is a Feller Markov semigroup. 
\beg{enumerate} \item[$(1)$] {\bf (D1)} implies $\{L^\nu_t\}_{\nu\in \scr P(E)}\in LDP_u(J)$, and the inverse holds when $E$ is locally compact. If moreover $P_t$ is strong Feller and $\mu$-irreducible   for some $t>0$, then $\{L^\nu_t\}_{\nu\in \scr P(E)}\in LDP(J)$ if and only if {\bf (D1)} holds.
 \item[$(2)$] {\bf (D2)} implies
 $\{L^x_t\}_{x\in D}\in LDP_u(J)$ for any compact set $D\subset E$, and the inverse holds provided   $E$ is locally compact. If $P_t$ is strong Feller and $\mu$-irreducible   for some $t>0$, then $\{L^x_t\}_{x\in D}\in LDP(J)$ for compact $D\subset E$ if and only if {\bf (D2)} holds.
\end{enumerate}
\end{thm}

Finally, we introduce the following approximation lemma which is easy to prove but useful in  applications, see for instance \cite[Theorems 4.2.16, 4.2.23]{DZ}, and see also \cite[Theorem 3.2]{RWW06} for a stronger version   called generalized contraction principle.

\beg{lem}[Approximation Lemma for LDP] \label{L1.1} Let $\{(L_t^\nu)_{t>0}, (\bar L_t^\nu)_{t>0}: \nu\in \scr I\}$    be two families of stochastic processes  on a Polish space $(E,\rr)$ for an index set $\scr I$.
If   $(\bar L^\nu_t)_{\nu\in scr I}\in LDP_u (J)($respectively $LDP_l(J)$)  and
 $$\lim_{t\to\infty}\ff 1 {\lll(t)} \sup_{\nu\in \scr A}\log \P(\rr(L_t^\nu,\bar L_t^\nu)>\dd)=-\infty,\ \ \dd>0,$$ then   $(L_t^\nu)_{\nu\in \scr I}\in LDP_u(J) ($respectively $LDP_l(J))$.\end{lem}

\section{Proofs of main results }

 To establish the LDP for $L_t^\nu$,  we will compare    \eqref{E'} with a reference equation:
\beq\label{BE} \d \bar X^\nu (t)= \big\{A\bar X^\nu (t) + \bar b(\bar X_t^\nu)\big\}\d t +\bar\si (\bar X_t^\nu) \d W(t),\ \ \bar X_0^\nu= X_0^\nu,\end{equation}
where $\bar b: \C\to\H,\ \bar\si: \C\to\mathbb L(\H)$ are measurable such that this equation has a unique mild segment solution for any initial value in $\C$, which is thus   a Markov process on $\C$. In applications, the coefficients in \eqref{BE} will be given by  the limit of $b_t(\cdot,\nu_t)$ and $\si_t(\cdot,\nu_t)$ as $t\to\infty$, where $b_t$ and $\si_t$ are in \eqref{E'} and $\nu_t:=\L_{X_t^\nu}$. Now, let
$$\bar L_t^\nu=\ff 1 t \int_0^t\dd_{\bar X_s^\nu}\d s,\ \ t>0.$$
We have the following result.

\beg{thm}\label{TM} Assume  that $\eqref{E'}$ and $\eqref{BE}$ are well-posed for any initial value $X_0$ with $\L_{X_0}\in \scr I$ and
$\L_{\bar X_0}\in \Psi(\scr I)$ respectively, where $\scr I$ is a non-empty subset of $\scr P(\C)$ and $\Psi:\scr  I\to\scr P(\C)$ is a map.   If $\{\bar L_t^\nu\}_{\nu\in \Psi(\scr I)}\in LDP_u(J) ($respectively $LDP_l(J))$ under the weak topology,   and
\beq\label{QQP} \sup_{\nu\in  \scr I}  \E\big[\e^{N\int_0^\infty \{\|X_s^\nu-\bar X_s^{\Psi(\nu)}\|_\infty\land 1\}\d s }\big]<\infty,\ \  N\ge 1,\end{equation}then    $\{L_t^\nu\}_{\nu\in \scr I}\in LDP_u(J)($respectively $LDP_l(J))$ under the weak topology.
 \end{thm}
\beg{proof} Consider the probability distance
\beq\label{RR} \rr(\LL_1,\LL_2):= \inf_{\Pi\in \C(\LL_1,\LL_2)} \int_{\C\times\C}\big\{\|\xi-\eta\|_\infty\land 1\big\}\Pi(\d\xi,\d\eta)\end{equation}
 on $\scr P(\C)$. It is well known that $\rr$ induces the weak topology on $\scr P(\C)$. Since
 $$\rr(L_t^\nu,\bar L_t^{\bar\nu})\le \ff 1 t \int_0^t \big\{\|X_s^\nu-\bar X_s^{\Psi(\nu)}\|_\infty\land 1\big\} \d s,\ \ t>0,$$
  \eqref{QQP} implies
 \beg{align*} &\lim_{t\to\infty}\ff 1 {t} \sup_{\nu\in \scr I}\log \P(\rr(L_t^\nu,\bar L_t^{\Psi(\nu)})>\dd)\\
 &\le \lim_{t\to\infty}\ff 1 {t} \sup_{\nu\in \scr I}\log \P\Big(N\int_0^t \big\{\|X_s^\nu-\bar X_s^{\Psi(\nu)}\|_\infty\land 1\big\} \d s> tN\dd\Big)\\
 &\le - N\dd,\ \ N\ge 1,\dd>0.\end{align*}Therefore,
 $$\lim_{t\to\infty}\ff 1 {t} \sup_{\nu\in \scr I}\log \P(\rr(L_t^\nu,\bar L_t^{\Psi(\nu)})>\dd)=-\infty,\ \ \dd>0.$$
  Then the desired  assertion follows from   Lemma \ref{L1.1} with   $\bar L^{\Psi(\nu)}_t$ replacing $\bar L_t^\nu$.  \end{proof}

\subsection{Proof of Theorem \ref{T01}}
Obviously, condition {\bf $(H_1)$} implies that the SDE
$$\d \bar X(t)= b(\bar X(t),\bar\mu)\d t +\si(\bar\mu)\d W(t)$$ is well-posed and the solution is a Markov Feller process, where $\bar\mu$ is the unique invariant probability measure of $P_t^*$. Let $\bar X^x(t)$ denote the solution starting at $x$.  According to Theorem \ref{T3}   and Theorem \ref{TM}, we only need to prove the following assertions:
 \beg{enumerate} \item[(a)]  For any $\lll>0$, there exist a constant $s>0$ and compact set $K\subset \R^d$, such that  \eqref{Wu2} holds for any compact set $K'\subset \R^d$ and
 $$\tau_K^x:=\inf\{t\ge 0: \bar X^x(t)\in K\},\ \ x\in \R^d.$$
\item[(b)] For any  $N\ge 1,$
$$ \sup_{\nu\in \scr B_{r,R} } \E\e^{N\int_0^\infty \{1\land |X^\nu(s)-\bar X^0(s)|^2\}\d s}<\infty.$$
\item[(c)] Under \eqref{ABC}, for any $\lll>0$ there exists a compact set $K\subset \R^d$ such that \eqref{Wu1} holds for $\bar X$, and
$$\sup_{\nu\in\scr P_2(\R^d)}   \E\e^{N\int_0^\infty \{1\land |X^\nu(s)-\bar X^\nu(s)|^2\}\d s}<\infty,\ \ N\ge 1.$$
\end{enumerate}
Indeed, by Theorem \ref{T3}(2), (a) implies the upper LDP (LDP if $\bar P_t$ is strong Feller and $\bar\mu$-irreducible) for $\bar L_t^x$ locally uniformly in $x$, in particular, $L_t^0$ satisfies the upper LDP (LDP if $\bar P_t$ is strong Feller and $\bar\mu$-irreducible).
Combining this with (b) and Theorem \ref{TM} for $\scr I=\B_{r,R}$ and $\Psi(\nu):=\dd_0,$   we prove the desired assertion   for $L_t^\nu$ with $\nu\in \scr B_{r,R}$. Finally, by Theorem \ref{T3}(1) and Theorem \ref{TM} with $\scr I=\scr P_2(\R^d)$ and $\Psi(\nu)=\nu$, (c) implies the upper LDP (LDP if $\bar P_t$ is strong Feller and $\bar\mu$-irreducible) for $L_t^\nu$ uniformly in $\nu\in \scr P_2(\R^d).$

\paragraph{Proof of (a).} By {\bf $(H_1)$}, there exist constants $\aa,\bb>0$ such that
\beq\label{*P1} \d |\bar X(t)|^2\le 2 \{\aa-\bb |\bar X(t)|^2\}\d t+ 2\<\bar X(t),\si(\bar\mu)\d W(t)\>.\end{equation}
Let   $\theta=\|\si\|^2_\infty$. Then  for any $\vv\in (0, \bb/\theta)$, there exist constants $c_1,c_2>0$ such that
\beg{align*} \d \e^{\vv|\bar X(t)|^2} &\le 2\vv\big\{\aa-(\bb-\vv \theta) |\bar X(t)|^2\big\}\e^{\vv|\bar X(t)|^2}\d t + \d M(t)\\
&\le \Big\{c_1-c_2\e^{\vv|\bar X(t)|^2}\big\}\d t +\d M(t) \end{align*} for some martingale $M(t)$. So,
\beq\label{*P2} \E\e^{\vv |\bar X^x(t)|^2} \le \e^{\vv |x|^2}+  \ff {c_1}{c_2},\ \ x\in\R^d.\end{equation}
To estimate $\tau_K^x$ for $K:=B_0(N)$, we take $N\ge N_0:=(2\aa/\bb)^{\ff 1 2}$. Then \eqref{*P1} implies
$$\d |\bar X(t)|^2\le  -\bb |\bar X(t)|^2 \d t+ 2\<\bar X^x(t),\si(\bar\mu)\d W(t)\>,\ \ t\le \tau_K^x.$$
For any $\dd>0$, we obtain
 \beg{align*} &\E\e^{\dd \int_0^{t\land\tau_K^x} |\bar X^x(s)|^2\d s}\le \e^{\dd\bb^{-1} |x|^2} \E\e^{2\dd\bb^{-1}\int_0^{t\land\tau_K^x} \<\bar X^x(s), \si(\bar X^x(s),\mu)\d W(s)\>}\\
&\le \e^{\dd \bb^{-1}|x|^2} \big(\E\e^{8\dd^2 \bb^{-2}\theta\int_0^{t\land \tau_K^x} |\bar X^x(s)|^2\d s}\big)^{\ff 12}.\end{align*}
Thus, taking $\dd\le \ff {\bb^2} {8 \theta}$ we arrive at
$$\E\e^{\dd N^2(t\land\tau_K^x)} \le \E\e^{\dd \int_0^{t\land\tau_K^x} |\bar X^x(s)|^2\d s}\le \e^{2\dd\bb^{-1} |x|^2}.$$
Letting $t\uparrow \infty$ implies
\beq\label{POO}  \E \e^{\dd N^2 \tau_K^x} \le \e^{2\dd\bb^{-1} |x|^2},\ \ x\in \R^d, N\ge N_0.\end{equation}
Combining this with   the Markov property and \eqref{*P2}, when $\dd\le \ff{\vv \bb}2$ we have
$$\E \e^{\dd N^2\tau_K^{\bar X^x(s)}}\le \E\e^{2\dd\bb^{-1} |\bar X^x(s)|^2}\le \e^{\vv |x|^2}+  \ff {c_1}{c_2},\ \ x\in\R^d, s\ge 0, N\ge N_0.$$
Therefore, for any $\lll>0$ there exists  compact $K\subset \R^d$ such that \eqref{Wu1}  holds.

\paragraph{Proof of (b).} Simply denote $X(t)=X^\nu(t), \bar X(t)= \bar X^0(t)$ and $\nu_t=\L_{X^\nu(t)}=P_t^*\nu$ for $\nu\in \B_{r,R}$. By {\bf $(H_1)$}, \eqref{EXPO} and It\^o's formula, we obtain
\beg{align*} \d |X(t)-\bar X(t)|^2 \le & \big\{-\kk_1|X(t)-\bar X(t)|^2 +\kk_2 \e^{-(\kk_1-\kk_2)t} \W_2(\bar\mu,\nu)^2\big\}\d t\\
&+ 2\<X(t)-\bar X(t), \{\si(\nu_t)-\si(\bar\mu)\}\d W(t)\>.\end{align*}
Letting $\gg(t)= \ff{|X(t)-\bar X(t)|^2 }{1+|X(t)-\bar X(t)|^2 }$, we derive
\beg{align*} &\d\log(1+|X(t)-\bar X(t)|^2)\le \big\{-\kk_1\gg(t)+\kk_2 \e^{-(\kk_1-\kk_2)t} \W_2(\bar\mu,\nu)^2\big\}\d t\\
&+ \ff{2}{1+|X(t)-\bar X(t)|^2} \<X(t)-\bar X(t), \{\si(\nu_t)-\si(\bar\mu)\}\d W(t)\>.\end{align*}
We deduce from this and \eqref{EXPO}  that for any $\lll>0,$
\beq\label{KK}\beg{split}  &\e^{-\ff{\lll \kk_2}{\kk_1-\kk_2} \W_2(\bar\mu,\nu)^2}\E\big[\e^{\lll\kk_1 \int_0^t \gg(s)\d s} \big]\\
& \le \E\Big[(1+|X_0|^2)^\lll\e^{\lll\int_0^t \ff{2}{1+|X(s)-\bar X(s)|^2}\<X(s)-\bar X(s), \{\si(\nu_s)-\si(\mu)\}\d W(s)\>}\Big]\\
&\le \E\bigg[(1+|X_0|^2)^\lll   \Big(\E\big[\e^{8\kk_2\lll^2  \int_0^t\gg (s) \W_2(\nu_s,\bar\mu)^2\d s} \big]\Big|\F_0\Big)^{\ff 1 2}\bigg]\\
&\le \big\{\nu\big((1+|\cdot|^2)^{2\lll}\big)\big\}^{\ff 1 2} \Big(\E\big[\e^{8\kk_2\lll^2  \W_2(\nu,\bar\mu)^2\int_0^t\gg (s) \e^{-(\kk_1-\kk_2)s}\d s }\big]\Big)^{\ff 1 2}\\
&\le C(\lll,R)  \Big(\E\big[\e^{\lll\kk_1  \int_0^t\gg (s)\d s } \big]\Big)^{\ff 1 2},\ \ t>0  \end{split}\end{equation}
holds for some constant $ C(\lll,R)>0$, where the last step is due to $\gg(s)\le 1$ and $\nu\in \B_{r,R}$. Therefore,
$$ \sup_{\nu\in \B_{r,R}} \E\Big[\e^{\lll\kk_1 \int_0^\infty \ff{|X^\nu(s)-\bar X^0(s)|^2 }{1+|X^\nu(s)-\bar X^0(s)|^2 }\d s }\Big] <\infty,\ \ \lll>0,$$
which implies   (b).

\paragraph{Proof of (c).}  Let \eqref{ABC} hold. Then there exist constants $c_1,c_2>0$ such that
\beq\label{PL} \d \e^{|\bar X(t)|^2}\le \big\{c_1 -c_2 |\bar X(t)|^{2+\vv}\e^{|\bar X(t)|^2}\big\} \d t+ 2\e^{|\bar X(t)|^2}\<\bar X(t),\si(\bar \mu)\d W(t)\>.\end{equation} This implies
  $$h_x(t):= \E \e^{|\bar X^x(t)|^2}\le c_1t +\e^{|x|^2} <\infty,\ \ t\ge 0, x\in \R^d. $$
Moreover, by  Jensen's inequality and the convexity of $[1,\infty)\ni r\mapsto r\log^{1+\vv/2}r,$ we deduce from \eqref{PL} that
$$h_x(t)\le h_x(0)+c_1t -c_2\int_0^t h_x(s) \log^{1+\vv/2}h_x(s)\d s,\ \ t\ge 0.$$
This and the comparison theorem imply   $h_x(t)\le \psi_x(t),$ where $\psi_x(t)$ solves the ODE
$$\psi'(t)= c_1-c_2\psi(t)\log^{1+\vv/2}\psi(t),\ \ \psi(0)=h_x(0)=\e^{|x|^2}.$$
So, \beq\label{GPP} \sup_{x\in\R^d} h_x(t)\le \sup_{\psi(0)\ge 1}\psi(t)=:c(t)<\infty.\end{equation}
On the other hand, by \eqref{PL},
 there exist constants $N_0,\bb>0$ such that for any $N\ge N_0$ and $K=B_0(N)$, we have
\beq\label{P*P} \d \e^{|\bar X^x(t)|^2} \le -\bb |\bar X^x(t)|^{2+\vv}\e^{|\bar X^x(t)|^2}\d t +2\e^{|\bar X^x(t)|^2}\<\bar X^x(t),\si(\bar\mu)\d W(t)\>,\ \ t\le \tau_K^x.\end{equation}
Combining this with \eqref{POO} and using the Markov property, when $2\dd\le\bb^2$ we arrive at
\beg{align*} &\E[\e^{\dd N^2\tau_K^x}]\le \e^{\dd N^2}+ \E\big[\e^{\dd N^2\tau_K^x}1_{\{\tau_K^x \ge 1\}}\big]\\
&\le \e^{\dd N^2} +\E\big[\e^{\dd N^2(1+\tau_K^{\bar X^x(1)})}1_{\{\tau_K^x \ge 1\}}\big]\\
&\le \e^{\dd N^2} (1+  \E\e^{|\bar X^x(1)|^2})\le \e^{\dd N^2}(1+c(1))<\infty,\ \ x\in\R^d, N\ge N_0.\end{align*}
Therefore, for any $\lll>0$, there exists compact set $K$ such that \eqref{Wu1}  holds.

Finally, repeating the proof of \eqref{GPP} using $X^\nu(t)$ replacing $\bar X^x(t)$, we derive
$$\sup_{\nu\in \scr P_2(\R^d)} \E [\e^{|X^\nu(1)|^2}]<\infty.$$ This together with \eqref{GPP} yields
\beq\label{GPPO} \sup_{\nu\in \scr P_2(\R^d)} \E \big[\e^{|X^\nu(1)|^2}+\e^{|\bar X^\nu(1)|^2}\big]<\infty.\end{equation}
On the other hand, as in \eqref{KK} but integrating from time $1$, we obtain
 \beg{align*}  &\e^{-\ff{\lll \kk_2}{\kk_1-\kk_2} \W_2(\bar\mu,\nu)^2}\E\big[\e^{\lll\kk_1 \int_1^t \ff{|X^\nu(s)-\bar X^\nu(s)|^2}{1+|X^\nu(s)-\bar X^\nu(s)|^2}\d s} \big]\\
& \le \E\Big[(1+|X^\nu(1)-\bar X^\nu(1)|^2)^\lll\e^{\lll\int_1^t \ff{2}{1+|X^\nu(s)-\bar X^\nu(s)|^2}\<X^\nu(s)-\bar X^\nu(s), \{\si(\nu_s)-\si(\bar \mu)\}\d W(s)\>}\Big]\\
&\le  \big\{ \E\big[(1+|X^\nu(1)-\bar X^\nu(1)|^2)^{2\lll} \big]  \big\}^{\ff 1 2}   \Big(\E\big[\e^{\lll\kk_1 \W_2(P_1^*\nu,\bar\mu)^2 \int_1^t\ff{|X^\nu(s)-\bar X^\nu(s)|^2\e^{-(\kk_1-\kk_2)s}}{1+|X^\nu(s)-\bar X^\nu(s)|^2}\d s } \big]\Big)^{\ff 1 2},\ \ t>1.  \end{align*}  
Combining this with \eqref{GPPO}, we   derive
$$ \sup_{\nu\in\scr P_2(\R^d)} \E\e^{\lll\kk_1 \int_1^\infty \ff{|X^\nu(s)-\bar X^\nu(s)|^2}{1+|X^\nu(s)-\bar X^\nu(s)|^2}\d s} <\infty,\ \   \lll\ge 1.$$
 Therefore, (c) holds.

\subsection{Proof of Theorem \ref{T02}}
As explained in the beginning of Subsection 4.1 that we only need to verify (a) and (b)  in the last subsection for the present model.
Comparing with the finite-dimensional case, the main difficulty is that  bounded sets are no longer compact.
To construct compact sets, let $\{e_i\}_{i\ge 1}$ be the eigenbasis of $A$; i.e. it is an orthonromal basis of $\H$  such that $Ae_i=-\lll_i e_i, i\ge 1.$  For any $N>0$, the set
$$K:=B_{0,\gg}(N)=\Big\{x\in \H: |x|_\gg^2:= \sum_{i=1}^\infty\<x,e_i\>^2\lll_i^\gg\le N^2\Big\}$$ is a compact set in $\H$.

\paragraph{Proof of (a).} Simply denote $\bar X(t)=\bar X^x(t)$  and $\tau_K=\tau_K^x:=\inf\{t\ge 0: \bar X^x(t)\in K\}$.  By  {\bf $(H_2)$} and \eqref{CCO}, we  may apply It\^o's formula to
  $$ \psi(\bar X(t)):=  \<(-A)^{\gg-1}\bar X(t), \bar X(t)\>=\sum_{i=1}^\infty \<\bar X(t),e_i\>^2 \lll_i^{\gg-1},$$
  such that for some constants $d_1,d_2>0$
  \beq\label{SSS} \d \psi(\bar X(t)) \le  (d_1-d_2 |\bar X(t)|_\gg^2)\d t + \d M(t),\end{equation}
  where $M(t):=2\sum_{i=1}^\infty \lll_i^{\gg-1}\<\bar X(t),e_i\> \<\si(\bar\mu)\d W(t), e_i\>$ for an orthonormal basis $\{e_i\}_{i\ge 1}$ of $\H$.
  Let $N\ge N_0:= (2d_1/d_2)^{\ff 1 2},$   and consider $\tau_K$ for $K= B_{0,\gg}(N)$.
  Then
  \beq\label{S**} d_1-d_2|\bar X(t)|_\gg^2 \le - d_1 |\bar X(t)|^2_\gg,\ \ t\le \tau_K.\end{equation}
  Since $\si$ is bounded, by $(H_2)$ there exists a constant $c>0$ such that
  $$\<M\>(t)\le c\int_0^t |\bar X(s)|^2\d s,\ \ t\ge 0.$$
  So, letting $\tau_n:=\inf\{t\ge 0: |\bar X(t)|\ge n\}$, we deduce form \eqref{SSS} and \eqref{S**} that
\beg{align*} &\E \e^{\int_0^{t\land \tau_n\land\tau_K}\dd d_1 |X(s)|_\gg^2\d s} \le \e^{\dd   \psi(x)} \big(\E\e^{2\dd^2\<M\>(t\land \tau_n\land\tau_K)}\big)^{\ff 1 2}\\
&\le  \e^{ \dd \psi(x)} \big(\E\e^{2c\dd^2\int_0^{t\land\tau_n\land\tau_K} |\bar X(s)|_\gg^2\d s} \big)^{\ff 1 2}<\infty,\  \ n\ge 1.\end{align*}
  Taking $\dd\le (2c)^{-1}$ leads to
$$\E\e^{\dd d_1N^2(t\land\tau_n\land\tau_K)} \le  \E \e^{\int_0^{t\land \tau_n\land\tau_K}\dd d_1  |\bar X(s)|_\gg^2\d s} \le \e^{ 2\dd \psi(x)},\ \ t\ge 0, n\ge 1.$$ Letting $t,n\to\infty $ we derive
$$\E\e^{\dd N^2 d_1  \tau_K} \le  \e^{ 2\dd \psi(x)},\ \ x\in \H.$$
  Combining this with the Markov property, we obtain
$$\E\e^{\dd N^2 d_1 \tau_K^{\bar X(s)}}\le \E  \e^{ 2\dd \psi(\bar X(s))},$$
and it is easy to see from \eqref{SSS} that the upper bound is locally bounded in $x$ when $\dd$ is small enough. Therefore, condition (a) is satisfied, since $N\ge N_0$ is arbitrary.

\paragraph{Proof of (b).}
By {\bf $(H_2)$}, \eqref{EXPP}, and It\^o's formula, we have
\beg{align*} \d |X^\nu(t)-\bar X^0(t)|^2\le &\big\{-2(\lll_1-\aa_1)|X^\nu(t)-\bar X^0(t)|^2+2\aa_2 \W_2(P_t^*\nu,\bar\mu)^2 \big\}\d t\\
& \ + 2\<X^\nu(t)-\bar X^0(t), \{\si(P_t^*\nu)-\si(\bar\mu)\}\d W(t)\>.\end{align*}The remainder of the proof is completely similar to that  of (b)  in the last subsection.

\subsection{Proof of Theorem \ref{TL1}}
Let $\theta\in [0,\lll_1]$ such that $\kk_p= \theta-(\aa_1+\aa_2)\e^{p\theta r_0}.$ 
\subsubsection{Proof of Theorem \ref{TL1}(1) }  For any $\nu_1,\nu_2\in \scr P_p(\C)$, take $X_0^{\nu_i} \in L^p(\OO\to\C,\F_0,\P)$
such that $\L_{X_0^{\nu_i}}=\nu_i,i=1,2,$ and
\beq\label{OP} \E\big[\|X_0^{\nu_1}-X_0^{\nu_2}\|_\infty^p\big]=\W_p(\nu_1,\nu_2)^p.\end{equation}
Since $\si$ is constant, we have
$$\d (X^{\nu_1}(t)-X^{\nu_2}(t))= \big\{A(X^{\nu_1}(t)-X^{\nu_2}(t)) + b(X_t^{\nu_1}, P_t^*\nu_1)-b(X_t^{\nu_2}, P_t^*\nu_2) \big\}\d t,\ \ t\ge 0.$$
By {\bf $(H_3)$} and noting that $\theta\in [0,\lll_1]$, we obtain
\beg{align*} &\d\big\{|X^{\nu_1}(t)-X^{\nu_2}(t)|^p\e^{p\theta t}\big\}\\
&= p\e^{p\theta t} |X^{\nu_1}(t)-X^{\nu_2}(t)|^{p-2}\big\{\<X^{\nu_1}(t)-X^{\nu_2}(t),(\theta+A)(X^\mu(t)-X^{\nu_2}(t))\>\\
 &\qquad\qquad\qquad \qquad\qquad \qquad + \<X^{\nu_1}(t)-X^{\nu_2}(t),b(X_t^{\nu_1}, P_t^*\nu_1)-b(X_t^{\nu_2}, P_t^*\nu_2)\>\big\}\d t\\
&\le p|X^{\nu_1}(t)-X^{\nu_2}(t)|^{p-1}\e^{p\theta t} \big\{\aa_1\|X_t^{\nu_1}-X_t^{\nu_1}\|_\infty + \aa_2\W_p(P_t^*\nu_1, P_t^*\nu_2)\big\}\d t,\ \ t\ge 0.\end{align*}
Letting $\psi(t)= \|X^{\nu_1}_t-X^{\nu_2}_t\|_\infty^p\e^{p\theta t}$, we derive
\beq\label{*DP1}\beg{split} &\psi(t) \le \e^{p\theta r_0} \sup_{s\in [(t-r_0),t]} |X^{\nu_1}(s)-X^{\nu_2}(s)|^p\e^{p\theta s} \\
&\le \e^{p\theta r_0} \|X_0^{\nu_1}-X_0^{\nu_2}\|_\infty^p +p\e^{p\theta r_0} \int_0^t   \big\{\aa_1\psi(s) + \aa_2\e^{\theta s} \W_p(P_s^*\nu_1, P_s^*\nu_2) \psi(s)^{\ff{p-1}p}\big\}\d s.\end{split}\end{equation} 
Combining this with \eqref{OP} and
$$\W_p(P_t^*\nu_1, P_t^*\nu_2)^p\le \E\|X_t^{\nu_1}-X_t^{\nu_2}\|_\infty^p,\ \ t\ge 0,$$
we arrive at
$$\E[\psi(t)]\le \e^{p\theta r_0} \W_p(\nu_1,\nu_2)^p + p\e^{p\theta r_0} (\aa_1 + \aa_2)\int_0^t   \E[\psi(s)] \d s,\ \ t\ge 0.$$
By Theorem \ref{EXU} we have $\E[\psi(t)]<\infty, t>0$. Then  Gronwall's lemma yields
$$\E[\psi(t)]\le  \{ \W_p(\nu_1,\nu_2)\}^p\e^{p\theta r_0+p(\aa_1+\aa_2)\e^{p\theta r_0} t} ,\ \ t\ge 0.$$
Therefore,
$$\W_p(P_t^*\nu_1, P_t^*\nu_2)^p\le \e^{-p\theta t}  \E[\psi(t)]\le \{\W_p(\nu_1,\nu_2)\}^p\e^{p\theta r_0-p\kk_p t} ,\ \ t\ge 0.$$

When $\kk_p>0$, it is standard that \eqref{ES1} implies the existence and uniqueness of   $P_t^*$-invariant probability measure $\mu$ such that
\eqref{ES2} holds, see, for instance,  \cite[Proof of Theorem 3.1(2)]{W18}.

\subsubsection{Proof of Theorem \ref{TL1}(2)} Let $\kk_p>0.$  To prove the LDP, let $\bar P_t$ be the Markov semigroup for the stationary equation \eqref{**} and   consider the LDP for $\bar L_t^\nu$. Since $\lll>0$ implies
$\sup_{r\in [0,\lll_1]} (r-\aa_1 \e^{rr_0})>0$ and noting that   {\bf $(H_3)$} implies
$$|b(\xi,\mu)-b(\eta,\mu)|\le \aa_1 \|\xi-\eta\|_\infty,$$
by \cite[Theorem  1.2]{BWY15} and $\kk_1\ge \kk_p>0$, the Markov semigroup $\bar P_t$ is hypercontractive. Thus, by the semigroup property and the interpolation theorem,
for any $q>1$ there exists $t_q>0$ such that $\bar P_{t_q} $ is uniformly integrable in $L^q(\mu)$. Moreover,
according to \cite[Theorem 4.2.4]{W13}, assumption {\bf $(H_3)$} implies that  for any $t>r_0$, there exists a constant
$c>0$ such that the following Harnack inequality holds:
\beq\label{HH1} \big(\bar P_{t_0} f(\eta)\big)^2\le (\bar P_{t_0} f^2(\xi)))
\e^{c\|\xi-\eta\|_\infty^2},\ \ \ \xi,\eta\in\C, f\in
\B_b(\C).\end{equation}
Obviously, $\bar\mu$ is also $\bar P_t$-invariant, then for any  $B\in\B(\C)$ such that $\bar\mu(B)>0$, we have $\bar\mu(\bar P_{t_0}1_B)=\mu(B)>0$, so that there exits $\eta\in\C$ such that $\bar P_{t_0}1_B(\eta)>0$. Then \eqref{HH1} implies $\bar P_{t_0}1_B(\xi)>0$ for all $\xi\in\C$, so that
$\bar\mu(1_A \bar P_{t_0} 1_B)>0$ for $\bar\mu(A),\bar\mu(B)>0,$ i.e. $\bar P_{t_0}$ is $\bar\mu$-irreducible. Therefore, by Theorem \ref{T2},
\beq\label{OUP} \bar L_t^\nu\in LDP(J)\ \text{ uniformly\ in\ }  \nu=h\bar\mu \in\scr P(\C) \ \text{with}\ \|h\|_{L^q(\bar\mu)}\le R,\ \ R>0. \end{equation}
Combining this with Theorem \ref{TM}, it remains to show that for any $\vv,R>0$, 
\beg{enumerate}
\item[(I)] $\{\bar L_t^\nu\}_{\nu\in \scr I_{\vv,R}} \in LDP(J);$
\item[(II)]  For any $\dd>0$, 
$$\lim_{t\to\infty} \ff 1 t \sup_{\nu\in \scr I_{\vv,R} }\log \P\bigg(\ff 1 t \int_0^t \{1\land \|X_s^\nu-\bar X_s^\nu\|_\infty\}\d s  >\dd\bigg)=-\infty.$$\end{enumerate}

\paragraph{For (I).}  Observing that for any $\xi,\eta\in \C$ we have 
$$\d (\bar X^\xi(t)-\bar X^\eta(t)) = \big\{A(\bar X^\xi(t)-\bar X^\eta(t)) + b(\bar X_t^\xi,\bar\mu)- b(\bar X^\eta_t,\bar\mu)\big\}\d t,$$
by the same reason leading to \eqref{*DP1} we obatin
$$\|\bar X_t^\xi-\bar X_t^\eta\|_\infty^p \e^{p\theta t}\le \e^{p\theta r_0} \|\xi-\eta\|_\infty^p +\aa_1p \e^{p\theta r_0} \|\bar X_s^\xi-\bar X_s^\eta\|_\infty^p \e^{p\theta s}\d s,\ \ t\ge 0.$$
Noting that $\kk_p\le \theta-\aa_1\e^{p\theta r_0},$ by Gronwall's inequality we get 
$$\|\bar X_t^\xi-\bar X_t^\eta\|_\infty^p\le \e^{p\theta r_0-p\{\theta- \aa_1\e^{p\theta r_0}\}t} \|\xi-\eta\|_\infty^p\le  
e^{p\theta r_0-p\kk_p\}t} \|\xi-\eta\|_\infty^p.$$
Combining this with \eqref{HH1} and using the semigroup property of $\bar P_t$, we find a constant $t_1>t_0$ such that 
$$\big(\bar P_{t_1} f(\eta)\big)^2\le (\bar P_{t_1} f^2(\xi)))
\e^{\vv\|\xi-\eta\|_\infty^2/2},\ \ \ \xi,\eta\in\C, f\in
\B_b(\C).$$ This implies that the invariant probability measure $\bar\mu$ has full support on $\C$, so that 
there exists a constant $c>0$ such that 
$$\sup_{\bar\mu(|f|^2)\le 1} (\bar P_{t_1} f(\xi))^2\le \ff 1 {\int_\C \e^{-\vv \|\xi-\eta\|_\infty^2/2}\bar\mu(\d\eta)} \le c \e^{\vv \|\xi\|_\infty^2},\ \ \xi\in \C.$$
Therefore, $\bar P_{t_1} $ has a density $p_{t_1}(\xi,\eta)$ with respect to $\bar \mu$ satisfying 
$$\int_\C p_{t_1}(\xi,\eta)^2 \bar\mu(\d\eta)\le c\e^{\vv\|\xi\|_\infty^2},\ \ \xi\in \C.$$
Consequently, for any $\nu\in \scr I_{\vv,R}$, $\bar\nu_{t_1}:=\L_{\bar X_{t_1}^\nu} $ has  density 
$$h(\eta):= \int_\C p_{t_1}(\xi,\eta)\nu(\d \xi)$$
with respect to $\bar\mu$ which  satisfies 
$$\bar\mu(|h|^2)\le \int_{\C\times \C}p_{t_1}(\xi,\eta)^2 \nu(\d\xi)\bar\mu(\d\eta)\le c \nu(\e^{\vv\|\cdot\|_\infty^2})\le cR.$$
Combining this with \eqref{OUP} and noting that  the Markov property of $\bar X_t$ implies that 
the law of $\bar L_t^{\bar\nu_{t_1}}$ coincides with that of 
$$\tt L_t^\nu:= \ff 1t\int_{t_1}^{t+t_1} \dd_{\bar X_s^\nu} \d s,$$
we prove 
 \beq\label{OUP'} \{\tt L_t^\nu\}_{\nu\in \scr I_{\vv,R}} \in LDP(J).  \end{equation}
On the other hand, for the distance $\rr$ in \eqref{RR} we have 
$$\rr(\tt L_t^\nu, \bar L_t^\nu)\le \ff {2t_1}t,\ \ t>0.$$ So, by Lemma \ref{L1.1} and \eqref{OUP'}  we prove (I).

\paragraph{For (II).}  By {\bf $(H_3)$}  and   \eqref{ES2}, there exist constants $c>0$ such that  for $\theta\in [0,\lll_1]$,
\beg{align*} \|X_t^\nu-\bar X^\nu_t\|_\infty\e^{\theta t} &\le \e^{\theta r_0} \sup_{s\in [(t-r_0)^+, t]} |X^\nu(s)-\bar X^\nu(s)|\e^{\theta s}\\
&\le \e^{\theta r_0} \int_0^t \e^{\theta s} \big\{\aa_1\|X_s^\nu-\bar X_s^\nu\|_\infty+\aa_2\W_p(P_s^*\nu, \mu)\big\}\d s\\
& \le c + \aa_1 \e^{\theta r_0} \int_0^t \e^{\theta s} \|X_s^\nu-\bar X_s^\nu\|_\infty\d s,\ \ t\ge 0, \nu\in\scr I_{\vv,R}.\end{align*}
By Gronwall's inequality we obtain
$$\sup_{\nu\in \scr I_{\vv,R}} \|X^\nu_t-\bar X^\nu_t\|_\infty\le c \exp\big[\{\aa_1 \e^{\theta r_0}- \theta\}t\big]\le c\e^{-\kk_p t},\ \ t>0.$$
This proves assertion (II).

\subsection{Proof of Theorem \ref{T03}}

By \eqref{K2}, we take $\theta\in (0,\lll_1]$ such that
\beq\label{K3}   \theta\e^{-\theta r_0}- K_2 -\aa' \|B\|> K_3.\end{equation}
For any $\aa>0$, let $$\rr_{\aa}(\xi_1,\xi_2):= \aa\|\xi_1^{(1)} -\xi_2^{(1)}\|_\infty + \|\xi_1^{(2)} -\xi_2^{(2)}\|_\infty,\  \ \xi_1,\xi_2\in\C.$$
We take  $X_0,Y_0\in L^2(\OO\to\C,\F_0,\P)$ such that $\L_{X_0}=\nu_1, \L_{Y_0}=\nu_2$ and
\beq\label{K4} \W_{p,\aa}(\nu_1,\nu_2)^p=\E\rr_\aa(X_0,Y_0)^p.\end{equation}

Let $X(t)$ and $Y(t)$ solves \eqref{E1} with initial values $X_0$ and $Y_0$ respectively. Then $(H_4^1)$ implies $A_1-\dd\le -\lll_1\le\theta$, so that
\beg{align*} &|X^{(1)}(t)-Y^{(1)}(t)|\le   |\e^{(A_1-\dd)t} \{X^{(1)}(0)-Y^{(1)}(0)\}|\\
&\qquad + \int_0^t \big|\e^{(A_1-\dd)(t-s)}\{\dd(X^{(1)}(s)-Y^{(1)}(s))+B(X^{(2)}(s)-Y^{(2)}(s))\}\big|\d s\\
&\le \e^{-\theta t}|X^{(1)}(0)-Y^{(1)}(0)|+\int_0^t\e^{-\theta(t-s)} \big\{\dd |X^{(1)}(s)-Y^{(1)}(s)|+\|B\|\cdot  |X^{(2)}(s)-Y^{(2)}(s)|\big\}\d s.\end{align*}
Equivalently,
\beg{align*} & \e^{\theta t} |X^{(1)}(t)-Y^{(1)}(t)|\le |X^{(1)}(0)-Y^{(1)}(0)| \\
&\qquad +\int_0^t \e^{\theta s} \big\{\dd |X^{(1)}(s)-Y^{(1)}(s)|+\|B\|\cdot  |X^{(2)}(s)-Y^{(2)}(s)|\big\}\d s.\end{align*}
Similarly, it follows from $A_2\le-\lll_1\le-\theta$ and $(H_4^2)$ that
\beg{align*}& \e^{\theta t} |X^{(2)}(t)-Y^{(2)}(t)|\le |X^{(2)}(0)-Y^{(2)}(0)| \\
&+\int_0^t \e^{\theta s} \big\{K_1 \|X^{(1)}_s-Y^{(1)}_s\|_\infty +K_2\|X^{(2)}_s-Y^{(2)}_s\|_\infty +K_3 \W_{p,\aa}(P_s^*\nu_1, P_s^*\nu_2)\big\}\d s.\end{align*}
Combining these with $\aa'\ge \aa$ and   that $\lll':= \ff 1 2\{\dd +K_2+\ss{(K_2-\dd)^2+4\|B\|}$ satisfies
$$\aa'\dd +K_1=\lll'\aa',\ \ \aa'\|B\|+K_2=\lll'>0,$$
we derive
\beg{align*} &\e^{\theta t} \rr_{\aa'}(X_t,Y_t)
 \le \e^{\theta r_0} \sup_{s\in [t-r_0,t]}\{\aa'|X^{(1)}(s)-Y^{(1)}(s)|+|X^{(2)}(s)-Y^{(2)}(s)|\}\e^{\theta s} \\
& \le \e^{\theta r_0} \rr_{\aa'}(X_0,Y_0)+ \e^{\theta r_0}\int_0^t \big\{(\dd\aa' +K_1)\|X_s^{(1)}-Y_s^{(1)}\|_\infty\\
&\qquad\qquad\qquad + (\aa'\|B\|+K_2)\|X_s^{(2)}-Y_s^{(2)}\|_\infty
  +K_3\W_{p,\aa}(P_s^*\nu_1,P_s^*\nu_2)\big\}\d s\\
&=\e^{\theta r_0} \rr_{\aa'}(X_0,Y_0)+ \e^{\theta r_0} \int_0^t \e^{\theta s} \big\{\lll'\rr_{\aa'}(X_s,Y_s)+K_3\E[\rr_{\aa}(X_s,Y_s)]\big\}\d s.\end{align*}
 By Gronwall's lemma, for $\kk:= \theta-\lll'\e^{\theta r_0} >0$ we have
\beg{align*}  \rr_{\aa'}(X_t,Y_t)  \le  \e^{\theta r_0-\kk t} \rr_{\aa'}(X_0,Y_0) + \e^{\theta r_0}K_3\int_0^t \e^{-\kk(t-s)}\E[\rr_{\aa}(X_s,Y_s)]\d s.\end{align*}
Therefore, for any $\vv>0$ there exists a constant $C(\vv)>0$ such that
$$\rr_{\aa'}(X_t,Y_t)^p\le C(\vv) \rr_{\aa'}(X_0,Y_0)^p \e^{-\kk p t} + \ff{K_3^p\e^{\theta r_0p}(1+\vv)}{\kk^{p-1}} \int_0^t \e^{-\kk(t-s)} \E[\rr_{\aa}(X_s,Y_s)^p]\d s.$$
Combining this with $\rr_{\aa'}\ge \rr_{\aa}$ and $\E[\rr_{\aa'}(X_t,Y_t)^p]<\infty$ due to Theorem \ref{EXU}, we deduce from this and Gronwall's lemma that
\beg{align*} & \W_{p,\aa}(P_t^*\nu_1,P_t^*\nu_2)^p \le \E [\rr_{\aa'}(X_t,Y_t)^p ]\\
 &\le  \ff{\aa^p C(\vv)}{(\aa')^p} \W_{p,\aa}(\nu_1,\nu_2)^p\exp\big[-\big(\kk -(1+\vv) K_3^p\e^{\theta r_0p} \kk^{1-p}\big)t\big].\end{align*}
It is easy to see that   \eqref{K3} implies $\kk > K_3^p\e^{\theta r_0p} \kk^{1-p}$, so that by taking small enough $\vv>0$ we prove
$$\W_p(P_t^*\nu_1,P_t^*\nu_2)\le c_1\e^{-c_2 t},\ \ t\ge 0, \nu_1,\nu_2\in \scr P_p(\C)$$ for some constants $c_1,c_2>0$. Consequently,    $P_t^*$ has a unique invariant probability measure $\bar\mu$  such that \eqref{K**} holds.

Similarly, by {\bf $(H_4)$} and \eqref{K**}, we find a constant $C>0$ such that for any   $X_0^\nu=\bar X_0^\nu\in L^p(\OO\to\C,\F_0,\P),$
  $$\int_0^\infty \|X_t^\nu-\bar X_t^\nu\|_\infty^2 \d t\le C,\ \ \nu\in\scr I_{R,q}.$$
  Moreover, it is easy to see that \eqref{K2}  implies the condition in \cite[Theorem 1.3]{BWY15} for the reference equation with $\bar\mu$ replacing the distribution of solution, so that $\bar  P_t$ is hypecontractive (hence uniformly integrable in $L^p(\bar\mu)$ for any $p>1$) for large $t>0$, and   the Harnack in \cite[Lemma 4.1]{BWY15} implies \eqref{HH1}. Then the desired  LDP can be proved in the same way as in the proof of Theorem \ref{TL1}.

\end{document}